\documentclass{amsart}
\usepackage[square,numbers]{natbib}

\RequirePackage{amsthm,amsmath}
\usepackage{amsfonts,epsfig,graphicx}

\newtheorem{thm}{Theorem}
\newtheorem{lem}[thm]{Lemma}

\theoremstyle{definition}

\newtheorem{ex}[thm]{Example}

\providecommand{\abs}[1]{\lvert#1\rvert}
\providecommand{\Abs}[1]{\Bigl\lvert#1\Bigr\rvert}
\providecommand{\norm}[1]{\lVert#1\rVert}
\newcommand\independent{\protect\mathpalette{\protect\independenT}{\perp}}
\def\independenT#1#2{\mathrel{\rlap{$#1#2$}\mkern2mu{#1#2}}}

\begin{document}
\title{Generating knockoffs via conditional independence}

\author{Emanuela Dreassi}
\address{Emanuela Dreassi, Dipartimento di Statistica, Informatica, Applicazioni (DiSIA), Universit\`a degli Studi di Firenze, Italy}
\email{emanuela.dreassi@unifi.it}

\author{Fabrizio Leisen}
\address{Fabrizio Leisen, University of Nottingham, UK}
\email{fabrizio.leisen@gmail.com}

\author{Luca Pratelli}
\address{Luca Pratelli, Accademia Navale, Livorno, Italy} \email{pratel@mail.dm.unipi.it}

\author{Pietro Rigo}
\address{Pietro Rigo, Dipartimento di Scienze Statistiche, Universit\`a di Bologna, Italy}
\email{pietro.rigo@unibo.it}



\begin{abstract}
Let $X$ be a $p$-variate random vector and $\widetilde{X}$ a knockoff copy of $X$ (in the sense of \cite{CFJL18}). A new approach for constructing $\widetilde{X}$ (henceforth, NA) has been introduced in \cite{JSPI}. NA has essentially three advantages: (i) To build $\widetilde{X}$ is straightforward; (ii) The joint distribution of $(X,\widetilde{X})$ can be written in closed form; (iii) $\widetilde{X}$ is often optimal under various criteria. However, for NA to apply, $X_1,\ldots,X_p$ should be conditionally independent given some random element $Z$. Our first result is that any probability measure $\mu$ on $\mathbb{R}^p$ can be approximated by a probability measure $\mu_0$ of the form
$$\mu_0\bigl(A_1\times\ldots\times A_p\bigr)=E\Bigl\{\prod_{i=1}^p P(X_i\in A_i\mid Z)\Bigr\}.$$
The approximation is in total variation distance when $\mu$ is absolutely continuous, and an explicit formula for $\mu_0$ is provided. If $X\sim\mu_0$, then $X_1,\ldots,X_p$ are conditionally independent. Hence, with a negligible error, one can assume $X\sim\mu_0$ and build $\widetilde{X}$ through NA. Our second result is a characterization of the knockoffs $\widetilde{X}$ obtained via NA. It is shown that $\widetilde{X}$ is of this type if and only if the pair $(X,\widetilde{X})$ can be extended to an infinite sequence so as to satisfy certain invariance conditions. The basic tool for proving this fact is de Finetti's theorem for partially exchangeable sequences. In addition to the quoted results, an explicit formula for the conditional distribution of $\widetilde{X}$ given $X$ is obtained in a few cases. In one of such cases, it is assumed $X_i\in\{0,1\}$ for all $i$.
\end{abstract}

\maketitle

\section{Introduction}\label{intro}

\noindent One of the main problems, both in statistics and machine learning, is to identify the explanatory variables which are to be discarded, for they don't have a meaningful effect on the response variable. To formalize, let $X_1,\ldots,X_p,Y$ be real random variables, where $Y$ is regarded as the response variable and $X_1,\ldots,X_p$ as the explanatory variables. A {\em Markov blanket} is a minimal subset $\mathcal{S}\subset\{1,\ldots,p\}$ such that
$$Y\independent (X_i:i\notin\mathcal{S})\mid (X_i:i\in\mathcal{S}).$$
Under mild conditions, a Markov blanket $\mathcal{S}$ exists, is unique, and $\{1,\ldots,p\}\setminus\mathcal{S}$ can be written as
$$\{1,\ldots,p\}\setminus\mathcal{S}=\bigl\{i:\quad Y\independent X_i\mid (X_1,\ldots,X_{i-1},X_{i+1},\ldots,X_p)\bigr\};$$
see e.g. \cite[p. 558]{CFJL18} and \cite[p. 8]{EDW}. The problem mentioned above is to identify $\mathcal{S}$.

\medskip

\noindent To any selection procedure concerned with this problem, we can associate the false discovery rate $E\Bigl(\frac{\abs{\Hat{\mathcal{S}}\setminus\mathcal{S}}}{\abs{\Hat{\mathcal{S}}}\vee 1}\Bigr)$, where $\Hat{\mathcal{S}}$ denotes the estimate of $\mathcal{S}$ provided by the procedure. As in the Neyman-Pearson theory, those selection procedures which take the false discovery rate under control worth special attention.

\medskip

\noindent One such procedure has been introduced by Barber and Candes; see \cite{BC15}, \cite{BCS20}, \cite{BCJW21}, \cite{CFJL18}, \cite{GGZ}, \cite{SSC19}. Let
$$X=(X_1,\ldots,X_p).$$
Roughly speaking, Barber and Candes' idea is to create an auxiliary vector
$$\widetilde{X}=(\widetilde{X}_1,\ldots,\widetilde{X}_p),$$
called a {\em knockoff copy of} $X$, which is able to capture the connections among $X_1,\ldots,X_p$. Once $\widetilde{X}$ is given, each $X_i$ is selected/discarded based on the comparison between it and $\widetilde{X}_i$. Intuitively, $\widetilde{X}_i$ plays the role of a control for $X_i$, and $X_i$ is selected if it appears to be considerably more associated with $Y$ than its knockoff copy $\widetilde{X}_i$. This procedure is a recent breakthrough as regards variable selection. In addition to take the false discovery rate under control, it has other merits. In particular, it works whatever the conditional distribution of $Y$ given $X$. More precisely, for the knockoff procedure to apply, {\em one must assign $\mathcal{L}(X)$ but is not forced to specify} $\mathcal{L}(Y\mid X)$. (Here and in the sequel, for any random elements $U$ and $V$, we denote by $\mathcal{L}(U)$ and $\mathcal{L}(U\mid V)$ the probability distribution of $U$ and the conditional distribution of $U$ given $V$, respectively).

\medskip

\noindent Let us make precise the conditions required to $\widetilde{X}$. For each $i\in\{1,\ldots,p\}$ and each point $x\in\mathbb{R}^{2p}$, define $f_i(x)\in\mathbb{R}^{2p}$ by swapping $x_i$ with $x_{p+i}$ and leaving all other coordinates of $x$ fixed. Then, $f_i:\mathbb{R}^{2p}\rightarrow\mathbb{R}^{2p}$ is a permutation. For instance, for $p=2$, one obtains $f_1(x)=(x_3,x_2,x_1,x_4)$ and $f_2(x)=(x_1,x_4,x_3,x_2)$. In this notation, $\widetilde{X}$ is a {\em knockoff copy of} $X$, or merely a {\em knockoff}, if
$$\textbf{(i)}\quad f_i(X,\widetilde{X})\sim (X,\widetilde{X})\text{ for each }i\in\{1,\ldots,p\}\quad\text{and}\quad\textbf{(ii)}\quad \widetilde{X}\independent Y\mid X.$$

\medskip

\noindent For the knockoff procedure to apply, one must select $\mathcal{L}(X)$ and construct $\widetilde{X}$. However, obtaining $\widetilde{X}$ is not easy. Condition (ii) does not create any problems, for it is automatically true whenever $\widetilde{X}$ is built based only on $X$, neglecting any information about $Y$. On the contrary, condition (i) is quite difficult to be realized. Current tractable methods to achieve (i) require conditions on $\mathcal{L}(X)$. To our knowledge, such methods are available only when $X$ is Gaussian \cite{CFJL18}, or the set of observed nodes in a hidden Markov model \cite{SSC19}, or conditionally independent given some random element \cite{JSPI} and \cite{GGZ}. The third condition (conditional independence) is discussed in Section \ref{n7x4r} and includes the other two as special cases. There are also some universal algorithms, such as the {\em Sequential Conditional Independent Pairs} \cite{CFJL18} and the {\em Metropolized Knockoff Sampler} \cite{BCJW21}, which are virtually able to cover any choice of $\mathcal{L}(X)$. However, these algorithms do not provide a closed formula for $\widetilde{X}$. More importantly, they are computationally intractable as soon as $\mathcal{L}(X)$ is complex; see \cite{BCJW21} and \cite{GGZ}. As a matter of fact, they work effectively only for some choices of $\mathcal{L}(X)$ (such us graphical models) but not for all. A last remark is that, even if one succeeds to build $\widetilde{X}$, the joint distribution of the pair $(X,\widetilde{X})$ could be unknown. This is a further shortcoming. In fact, after observing $X=x$, it would be natural to sample a value $\widetilde{x}$ for $\widetilde{X}$ from the conditional distribution $\mathcal{L}(\widetilde{X}\mid X=x)$. But this is impossible if $\mathcal{L}(\widetilde{X}\mid X=x)$ is unknown.

\medskip

\noindent In a nutshell, the above remarks may be summarized as follows. {\em If $X$ is not conditionally independent} (in the sense of Section \ref{n7x4r}), then:

\medskip

\begin{itemize}

\item How to build a reasonable knockoff $\widetilde{X}$ is unknown.

\item The existing numerical algorithms are computationally heavy and may fail to work.

\item Even if one succeeds to build $\widetilde{X}$, the joint distribution of the pair $(X,\widetilde{X})$ is unknown.

\end{itemize}

\subsection{A new approach to knockoffs construction}\label{n7x4r} As noted above, while powerful and effective, the knockoff procedure suffers from some shortcomings due to the difficulty of building a reasonable knockoff $\widetilde{X}$. Such shortcomings are partially overcome by a new method for constructing $\widetilde{X}$, based on conditional independence, introduced in \cite{JSPI}. Similar ideas were also previously developed in \cite{GGZ}. Another related reference is \cite{BARJAN}. In this section, we recall the main features of this method.

\medskip

\noindent Suppose that $X_1,\ldots,X_p$ are conditionally independent given some random element $Z$. Denote by $\Theta$ the set where $Z$ takes values and by $\gamma$ the probability distribution of $Z$. Moreover, let $\mathcal{B}$ be the Borel $\sigma$-field on $\mathbb{R}$ and
$$P_i(A\mid\theta)=P(X_i\in A\mid Z=\theta)\quad\text{for all }i=1,\ldots,p,\,\theta\in\Theta\text{ and }A\in\mathcal{B}.$$
Note that $\gamma$ is a probability measure on $\Theta$ and each $P_i(\cdot\mid\theta)$ is a probability measure on $\mathbb{R}$. Since $X_1,\ldots,X_p$ are conditionally independent given $Z$,
\begin{gather}\label{ci88}
P(X_1\in A_1,\ldots,X_p\in A_p)=E\Bigl\{\prod_{i=1}^pP(X_i\in A_i\mid Z)\Bigr\}
\\=\int_\Theta\prod_{i=1}^p\,P_i(A_i\mid\theta)\,\gamma(d\theta)\quad\quad\text{for all }A_1,\ldots,A_p\in\mathcal{B}.\notag
\end{gather}
Hence, one can define a probability measure $\lambda$ on $\mathbb{R}^{2p}$ as
$$\lambda\bigl(A_1\times\ldots\times A_p\times B_1\times\ldots\times B_p\bigr)=\int_\Theta\,\prod_{i=1}^p\,P_i(A_i\mid\theta)\,P_i(B_i\mid\theta)\,\gamma(d\theta)$$
where $A_i\in\mathcal{B}$ and $B_i\in\mathcal{B}$ for all $i$. In \cite[Th. 12]{JSPI}, it is shown that any $p$-variate random vector $\widetilde{X}$ such that
$$\mathcal{L}(X,\widetilde{X})=\lambda$$
is a knockoff copy of $X$.

\medskip

\noindent Thus, arguing as above, not only one builds $\widetilde{X}$ in a straightforward way but also obtains the joint distribution of $(X,\widetilde{X})$, namely
\begin{gather}\label{g4rx3e7}
P(X_1\in A_1,\ldots,X_p\in A_p,\widetilde{X}_1\in B_1,\ldots,\widetilde{X}_p\in B_p)=
\\=\int_\Theta\prod_{i=1}^p\,P_i(A_i\mid\theta)\,P_i(B_i\mid\theta)\,\gamma(d\theta).\notag
\end{gather}
The price to be paid is to assign $\mathcal{L}(X)$ so as to satisfy \eqref{ci88}. (Recall that the choice of $\mathcal{L}(X)$ is a statistician's task). But this price is not expensive for two reasons. The first one is quite practical. The probability measures satisfying \eqref{ci88} are flexible enough to cover most real situations. Modeling $X_1,\ldots,X_p$ as conditionally independent (given some $Z$) is actually reasonable in a number of practical problems. The second reason is theoretical and is based on the results of this paper. Indeed, even if \eqref{ci88} fails, $\mathcal{L}(X)$ can be approximated arbitrarily well by probability measures satisfying \eqref{ci88}; see Theorems \ref{t1} and \ref{t2} below.

\medskip

\noindent The previous approach has two further advantages. First, $\widetilde{X}$ is often optimal under some criterions, such as mean absolute correlation and reconstructability. This is discussed in Example \ref{v63e0qa2}. However, we note by now that
$$\text{cov}(X_i,\widetilde{X}_i)=0\quad\text{if }Z\text{ is such that }E(X_i\mid Z)=0.$$
Second, even if it is not Bayesian from the conceptual point of view, the previous approach largely exploits Bayesian tools. Hence, to  construct $\widetilde{X}$ and evaluate $\mathcal{L}(X,\widetilde{X})$, all the Bayesian machinery can be recovered.

\medskip

\noindent To illustrate, suppose that $P_i(\cdot\mid\theta)$ admits a density $f_i(\cdot\mid\theta)$ with respect to some dominating measure $\lambda_i$. For instance, $\lambda_i$ could be Lebesgue measure or counting measure. Then, the joint densities of $X$ and $(X,\widetilde{X})$ are, respectively,
\begin{gather*}
h(x)=h(x_1,\ldots,x_p)=\int_\Theta\prod_{i=1}^p\,f_i(x_i\mid\theta)\,\gamma(d\theta)\quad\quad\text{and}
\\f(x,\widetilde{x})=f(x_1,\ldots,x_p,\widetilde{x}_1,\ldots,\widetilde{x}_p)
=\int_\Theta\prod_{i=1}^p\,f_i(x_i\mid\theta)\,f_i(\widetilde{x}_i\mid\theta)\,\gamma(d\theta)
\end{gather*}
where $x$ and $\widetilde{x}$ denote points of $\mathbb{R}^p$. In turn, assuming $h(x)>0$ for the sake of simplicity, the conditional density of $\widetilde{X}$ given $X=x$ can be written as
$$\frac{f(x,\widetilde{x})}{h(x)}=\frac{\int_\Theta\prod_{i=1}^p\,f_i(x_i\mid\theta)\,f_i(\widetilde{x}_i\mid\theta)\,\gamma(d\theta)}{\int_\Theta\prod_{i=1}^p\,f_i(x_i\mid\theta)\,\gamma(d\theta)}.$$
Therefore, we have an explicit formula for $\mathcal{L}(\widetilde{X}\mid X=x)$.

\medskip

\noindent In the rest of this paper, to make the exposition easier, a knockoff $\widetilde{X}$ obtained as above (i.e., a knockoff $\widetilde{X}$ satisfying equation \eqref{g4rx3e7}) is said to be a {\em conditional independence knockoff} (CIK). To highlight the connection between $\widetilde{X}$ and $X$, {\em we  also say that $\widetilde{X}$ is the CIK of $X$}.

\subsection{Content of this paper}\label{6b8m3d0k} This paper is basically a follow up of \cite{JSPI}. It consists of two results, two examples, and a numerical experiment. The results are of the theoretical type. They aim to characterize the CIKs, to show that they can be applied to virtually any real situation, and to highlight some of their optimality properties. The examples provide an explicit formula for $\mathcal{L}(\widetilde{X}\mid X=x)$ in two (meaningful) cases: mixtures of 2-valued (or 3-valued) distributions and mixtures of centered normal distributions. In particular, the first example deals with the case $X_i\in\{0,1\}$ for all $i$. Such a case is important in applications, mainly in a genetic framework. Nevertheless, apart from our example, we are not aware of any theoretical investigation of this case. Finally, in the numerical experiment, the CIKs are tested against simulated and real data.

\medskip

\noindent In the sequel, for any $d\ge 1$, a probability measure on $\mathbb{R}^d$ is called {\em absolutely continuous} if it admits a density with respect to Lebesgue measure on $\mathbb{R}^d$. Moreover, $\mathcal{P}$ is the class of all probability measures on $\mathbb{R}^p$ and $\mathcal{P}_0\subset\mathcal{P}$ is the subclass consisting of those $\mu_0\in\mathcal{P}$ of the form
$$\mu_0(A_1\times\ldots\times A_p)=\int_\Theta\prod_{i=1}^p\,P_i(A_i\mid\theta)\,\gamma(d\theta),$$
for some choice of $\Theta$, $\gamma$ and $P_i(\cdot\mid\theta)$ such that $P_i(\cdot\mid\theta)$ is absolutely continuous for all $i$ and $\theta$.

\medskip

\noindent We next briefly describe our two results. Moreover, by means of an example, we point out some optimality properties of the CIKs.

\medskip

\noindent Our first result (henceforth, R1) is that, for all $\mu\in\mathcal{P}$ and $\epsilon>0$, there is $\mu_0\in\mathcal{P}_0$ such that
\begin{gather*}
d_{BL}(\mu,\mu_0)<\epsilon\quad\text{and}\quad d_{TV}(\mu,\mu_0)<\epsilon\text{ if }\mu\text{ is absolutely continuous.}
\end{gather*}
In addition, an explicit formula for $\mu_0$ is provided. Here, $d_{BL}$ and $d_{TV}$ are the bounded Lipschitz metric and the total variation metric, respectively. Their definitions are recalled in Section \ref{cho99j5}.

\medskip

\noindent The motivation for R1 is that, to build a CIK, one needs $\mathcal{L}(X)\in\mathcal{P}_0$. This is not guaranteed, however, since the choice of $\mathcal{L}(X)$ is not subjected to any constraint. Hence, it is natural to investigate whether $\mathcal{L}(X)$ can be at least approximated by elements of $\mathcal{P}_0$. Because of R1, this is actually true. Roughly speaking, R1 aims to support $\mathcal{P}_0$ by showing that its elements are (approximatively) able to model any real situation.

\medskip

\noindent In addition to the previous motivation, R1 has also some practical utility. Suppose $\mathcal{L}(X)=\mu$ for some $\mu\in\mathcal{P}$. To fix ideas, suppose $\mu$ is absoutely continuous. If $\mu$ is arbitrary, how to build a reasonable knockoff $\widetilde{X}$ is unknown. However, given $\epsilon>0$, there is $\mu_0\in\mathcal{P}_0$ such that $d_{TV}(\mu,\mu_0)<\epsilon$. Such a $\mu_0$ can be built explicitly (recall that R1 provides an explicit formula for $\mu_0$). Denote by $T$ a $p$-variate random vector such that $\mathcal{L}(T)=\mu_0$. Since $\mu_0\in\mathcal{P}_0$, the CIK $\widetilde{T}$ of $T$ can be obtained straightforwardly. Then,
\begin{gather*}
d_{TV}\left(\mathcal{L}(\widetilde{X}),\,\mathcal{L}(\widetilde{T})\right)=d_{TV}(\mu,\mu_0)<\epsilon
\end{gather*}
for {\em any} knockoff copy $\widetilde{X}$ of $X$. Hence, by the robustness properties of the knockoff procedure \cite{BCS20}, $\widetilde{T}$ should be a reasonable approximation of $\widetilde{X}$.

\medskip

\noindent Our second result (henceforth, R2) is a characterization of the CIKs. Let $\mathcal{K}$ denote the class of the CIKs, that is
$$\mathcal{K}=\bigl\{\widetilde{X}:\,\mathcal{L}(X,\widetilde{X})\text{ admits representation \eqref{g4rx3e7} for some }\Theta,\,\gamma\text{ and }P_i(\cdot\mid\theta)\bigr\}.$$
Moreover, for any knockoff $\widetilde{X}$, say that $(X,\widetilde{X})$ is {\em infinitely extendable} if there is an (infinite) sequence $V=(V_1,V_2,\ldots)$ such that

\begin{itemize}

\item $(V_1,\ldots,V_{2p})\sim (X,\widetilde{X})$;

\item $V$ satisfies the same invariance condition as $(X,\widetilde{X})$ (this condition is formalized in Section \ref{e3hhu7}).
\end{itemize}

\noindent Then, R2 states that
\begin{gather*}
\widetilde{X}\in\mathcal{K}\quad\Leftrightarrow\quad (X,\widetilde{X})\text{ is infinitely extendable.}
\end{gather*}
Hence, if $(X,\widetilde{X})$ is required to be infinitely extendable, then $X$ {\em must} be conditionally independent (given some $Z$) and $\widetilde{X}$ {\em must} be the CIK of $X$. The proof of R2 is based on de Finetti's theorem for partially exchangeable sequences.

\medskip

\noindent Based on R2, a question is whether infinite extendability of $(X,\widetilde{X})$ is a reasonable condition. To answer, two facts are to be stressed. Firstly, by Finetti's theorem, infinite extendability of $(X,\widetilde{X})$ essentially amounts to conditional independence of $X$ and $\widetilde{X}$. Secondly, for the knockoff procedure to have a low type II error rate, it is desirable that $X$ and $\widetilde{X}$ are "as independent as possible"; see e.g. \cite[p. 563]{CFJL18} and \cite{SPJA}. Now, to have $X$ and $\widetilde{X}$ as independent as possible, a reasonable strategy is to take $X$ and $\widetilde{X}$ conditionally independent, or equivalently to require $(X,\widetilde{X})$ to be infinitely extendable.

\begin{ex}\label{v63e0qa2}\textbf{(Optimality of the CIKs).}
Suppose $E(X_i^2)<\infty$ and var$(X_i)>0$ for all $i$. Obviously, $\widetilde{X}$ should be selected so as to make the power of the knockoff procedure as high as possible. To this end, two criterions are to minimize the {\em mean absolute correlation}
$$\sum_{i=1}^p\,\Abs{\frac{\text{cov}(X_i,\widetilde{X}_i)}{\text{var}(X_i)}},$$
and to minimize the {\em reconstructability index}
\begin{gather*}
\sum_{i=1}^pE\bigl\{\text{var}(X_i\mid L_i)\bigr\}^{-1}\quad\,\text{where }\,L_i=(X_1,\ldots,X_{i-1},X_{i+1},\ldots,X_p,\widetilde{X}_1,\ldots,\widetilde{X}_p).
\end{gather*}
The first criterion (mean absolute correlation) is quite popular in the machine learning comunity. At least in some cases, however, it is overcome by the second (reconstructability index); see \cite{BCJW21} and \cite{SPJA}. Note also that
\begin{gather*}
E\bigl\{\text{var}(X_i\mid L_i)\bigr\}=E\bigl\{E(X_i^2\mid L_i)-E(X_i\mid L_i)^2\bigr\}=E(X_i^2)-E\bigl\{E(X_i\mid L_i)^2\bigr\}\le E(X_i^2).
\end{gather*}
Suppose now that $X_1,\ldots,X_p$ are conditionally independent, given some random element $Z$, and $\widetilde{X}$ is the CIK of $X$. Suppose also that $E(X_i\mid Z)=0$ a.s. for all $i$. Then,
$$\text{cov}(X_i,\widetilde{X}_i)=E(X_i\,\widetilde{X}_i)=E\bigl\{E(X_i\,\widetilde{X}_i\mid Z)\bigr\}=E\bigl\{E(X_i\mid Z)\,E(\widetilde{X}_i\mid Z)\bigr\}=0.$$
Therefore, $\widetilde{X}$ is optimal under the first criterion. Moreover,
$$E(X_i\mid L_i)=E\bigl\{E(X_i\mid Z,\,L_i)\mid L_i\bigr\}=E\bigl\{E(X_i\mid Z)\mid L_i\bigr\}=0\quad\quad\text{a.s.}$$
Hence, $E\bigl\{\text{var}(X_i\mid L_i)\bigr\}=E(X_i^2)$ and the reconstructability index attains its minimum value $\sum_{i=1}^pE(X_i^2)^{-1}$. Therefore, $\widetilde{X}$ is optimal under the second criterion as well.
\end{ex}

\section{Theoretical results}\label{cho99j5}

\noindent We first recall some (well known) definitions. A function $f:\mathbb{R}^p\rightarrow\mathbb{R}$ is said to be {\em Lipschitz} if there is a constant $b\ge 0$ such that
\begin{gather*}
\abs{f(x)-f(y)}\le b\,\norm{x-y}\quad\quad\text{for all }x,\,y\in\mathbb{R}^p,
\end{gather*}
where $\norm{\cdot}$ is the Euclidean norm. In this case, we also say that $f$ is $b$-Lipschitz or that $b$ is a Lipschitz constant for $f$.

\medskip

\noindent We remind that $\mathcal{P}$ denotes the class of all probability measures on $\mathbb{R}^p$. The {\em bounded Lipschitz metric} $d_{BL}$ and the {\em total variation metric} $d_{TV}$ are two distances on $\mathcal{P}$. If $\mu,\,\nu\in\mathcal{P}$, they are defined as
$$d_{BL}(\mu,\nu)=\sup_g\,\Abs{\int_{\mathbb{R}^p} g\,d\mu-\int_{\mathbb{R}^p} g\,d\nu}\quad\text{and}\quad d_{TV}(\mu,\nu)=\sup_A\,\abs{\mu(A)-\nu(A)}$$
where $\sup_g$ is over the 1-Lipschitz functions $g:\mathbb{R}^p\rightarrow [-1,1]$ and $\sup_A$ is over the Borel subsets $A\subset\mathbb{R}^p$. Among other things, $d_{BL}$ has the property that
$$\mu_n\rightarrow\mu\text{ weakly}\quad\Leftrightarrow\quad d_{BL}(\mu_n,\mu)\rightarrow 0$$
where $\mu_n,\,\mu\in\mathcal{P}$. We also note that $d_{BL}$ and $d_{TV}$ are connected through the inequality $d_{BL}\le 2\,d_{TV}$.

\medskip

\noindent We next turn to our main results.

\medskip

\subsection{$\mathcal{P}_0$ is dense in $\mathcal{P}$}\label{k9d555r} Let $\mathcal{P}_0$ be the class of those probability measures $\mu_0\in\mathcal{P}$ which can be written as
$$\mu_0(A_1\times\ldots\times A_p)=\int_\Theta\prod_{i=1}^p\,P_i(A_i\mid\theta)\,\gamma(d\theta),\quad\quad A_1,\ldots,A_p\in\mathcal{B},$$
for some choice of $\Theta$, $\gamma$ and $P_i(\cdot\mid\theta)$. To avoid trivialities, $P_i(\cdot\mid\theta)$ is assumed to be absolutely continuous for all $i=1,\ldots,p$ and $\theta\in\Theta$. The latter assumption is motivated by the next example.

\begin{ex}\textbf{(Why $P_i(\cdot\mid\theta)$ absolutely continuous).}
Suppose
\begin{gather}\label{z3ok8}
\Theta=\mathbb{R}^p,\,\gamma=\mathcal{L}(X)\text{ and }P_i(A\mid\theta)=1_A(\theta_i)
\end{gather}
where $\theta=(\theta_1,\ldots,\theta_p)$ denotes a point of $\mathbb{R}^p$. Then,
$$P(X_1\in A_1,\ldots,X_p\in A_p)=\int_{\mathbb{R}^p}\prod_{i=1}^p1_{A_i}(\theta_i)\,\gamma(d\theta)=\int_\Theta\prod_{i=1}^p\,P_i(A_i\mid\theta)\,\gamma(d\theta).$$
Hence, without some further constraint (such as $P_i(\cdot\mid\theta)$ absolutely continuous), one would obtain $\mathcal{P}_0=\mathcal{P}$ with $\Theta$, $\gamma$ and $P_i(\cdot\mid\theta)$ as in \eqref{z3ok8}. However, this is not practically useful. In fact, under \eqref{z3ok8}, the CIK $\widetilde{X}$ of $X$ is the trivial knockoff $\widetilde{X}=X$, which is unsuitable to perform the knockoff procedure.
\end{ex}

\medskip

\noindent If $\mathcal{L}(X)\in\mathcal{P}_0$, it is straightforward to obtain the CIK $\widetilde{X}$ of $X$ and to write $\mathcal{L}(X,\widetilde{X})$ in closed form. But clearly it may be that $\mathcal{L}(X)\notin\mathcal{P}_0$. In this case, it is quite natural to investigate whether $\mathcal{L}(X)$ can be approximated by elements of $\mathcal{P}_0$. This is actually possible and the approximation is very strong if $\mathcal{L}(X)$ is absolutely continuous.

\begin{thm}\label{t1}
For all $\mu\in\mathcal{P}$ and $\epsilon>0$, there is $\mu_0\in\mathcal{P}_0$ such that $d_{BL}(\mu,\mu_0)<\epsilon$. In particular, one such $\mu_0$ is
\begin{gather}\label{x48uh6tqa3}
\mu_0(A)=\int_{\mathbb{R}^p}\mathcal{N}_p(x,\,c I)(A)\,\mu(dx)\quad\quad\text{for all Borel sets }A\subset\mathbb{R}^p
\end{gather}
where $c=\epsilon^2/2p$ and $\mathcal{N}_p(x,\,c I)$ denotes the Gaussian law on $\mathbb{R}^p$ with mean $x $ and covariance matrix $c\,I$, i.e.
$$\mathcal{N}_p(x,\,c I)(A)=(2\,\pi\,c)^{-p/2}\int_A\exp\Bigl\{-\frac{\norm{y-x}^2}{2c}\Bigr\}\,dy.$$
\end{thm}

\medskip

\begin{thm}\label{t2}
Suppose $\mu\in\mathcal{P}$ is absolutely continuous. Then, for each $\epsilon>0$, there is $\mu_0\in\mathcal{P}_0$ such that $d_{TV}(\mu,\mu_0)<\epsilon$. Moreover, if $\mu$ has a Lipschitz density, one such $\mu_0$ can be defined by \eqref{x48uh6tqa3} with
$$c=\frac{1}{4p}\,\Bigl(\frac{\epsilon}{b\,m(B)}\Bigr)^2,$$
where $b$ is a Lipschitz constant for the density of $\mu$, $m$ is the Lebesgue measure on $\mathbb{R}^p$ and $B\subset\mathbb{R}^p$ is any Borel set satisfying $\mu(B^c)<\epsilon/2$ and $0<m(B)<\infty$.
\end{thm}

\medskip

\noindent Theorems \ref{t1} and \ref{t2} are proved in the Supplementary Material.

\medskip

\noindent It is worth noting that, in addition to \eqref{x48uh6tqa3}, there are other laws $\mu_0\in\mathcal{P}_0$  satisfying the inequalities $d_{BL}(\mu,\mu_0)<\epsilon$ or $d_{TV}(\mu,\mu_0)<\epsilon$. Moreover, in the second part of Theorem \ref{t2}, the Lipschitz condition on the density of $\mu$ can be weakened at the price of making $\mu_0$ slightly more involved.

\medskip

\noindent The motivation of Theorems \ref{t1}-\ref{t2} has been mentioned in Section \ref{6b8m3d0k}. In short, if $\mathcal{L}(X)\notin\mathcal{P}_0$, the CIK $\widetilde{X}$ of $X$ cannot be built. However, Theorems \ref{t1}-\ref{t2} imply that $\mathcal{L}(X)$ can be approximated by elements of $\mathcal{P}_0$. Hence, with a negligible error, it can be assumed $X\sim\mu_0$ and the CIK $\widetilde{X}$ of $X$ can be easily obtained. This is our main motivation. However, Theorems \ref{t1}-\ref{t2} have a practical implication as well. Suppose $\mathcal{L}(X)=\mu$ for some $\mu\in\mathcal{P}$. To fix ideas, suppose $\mu$ is absolutely continuous with a Lipschitz density. Fix $\epsilon>0$, define $\mu_0\in\mathcal{P}_0$ as in Theorem \ref{t2}, and call $T$ a $p$-variate vector such that $\mathcal{L}(T)=\mu_0$. Since $\mu_0\in\mathcal{P}_0$, the CIK $\widetilde{T}$ of $T$ can be easily built. Moreover, given {\em any} knockoff copy $\widetilde{X}$ of $X$, since $\widetilde{X}\sim X\sim\mu$ and $\widetilde{T}\sim T\sim \mu_0$, Theorem \ref{t2} yields
$$d_{TV}\left(\mathcal{L}(\widetilde{X}),\,\mathcal{L}(\widetilde{T})\right)=d_{TV}(\mu,\mu_0)<\epsilon.$$
Therefore, by the robustness properties of the knockoff procedure \cite{BCS20}, $\widetilde{T}$ is expected to be a reasonable approximation of $\widetilde{X}$. (Obviously, the latter claim should be supported by a numerical comparison of the power and the false discovery rate corresponding to $\widetilde{X}$ and $\widetilde{T}$. Such a comparison is not trivial, however, {\em since $\widetilde{X}$ is unknown for arbitrary} $\mu$).

\medskip

\noindent Finally, two remarks are in order. The first is summarized by the following lemma.

\begin{lem}\label{y7x4ww1}
Let $\widetilde{T}$ be the CIK of $T$, where $T\sim\mu_0$ with $\mu_0$ given by \eqref{x48uh6tqa3}. Then,
$$(T,\widetilde{T})\sim (L+M,\,L+N)$$
where $L,\,M,\,N$ are independent, $L\sim\mu$ and $M\sim N\sim\mathcal{N}_p(0,\,c I)$.
\end{lem}
\begin{proof} For any Borel sets $A,\,B\subset\mathbb{R}^p$, one obtains
\begin{gather*}
P(L+M\in A,\,L+N\in B)=\int_{\mathbb{R}^p}P(x+M\in A,\,x+N\in B)\,\mu(dx)
\\=\int_{\mathbb{R}^p}P(x+M\in A)\,P(x+N\in B)\,\mu(dx)
\\=\int_{\mathbb{R}^p}\mathcal{N}_p(x,\,c I)(A)\,\mathcal{N}_p(x,\,c I)(B)\,\mu(dx)=P(T\in A,\,\widetilde{T}\in B).
\end{gather*}
\end{proof}
\noindent Lemma \ref{y7x4ww1} makes clear the structure of $\mathcal{L}(T,\widetilde{T})$ and may be useful for sampling from such distribution.

\medskip

\noindent The second remark is that, if $\mathcal{L}(X,\widetilde{X})$ is absolutely continuous and has a Lipschitz density, the pair $(T,\widetilde{T})$ can be taken such that
$$d_{TV}\left(\mathcal{L}(X,\widetilde{X}),\,\mathcal{L}(T,\widetilde{T})\right)<\epsilon.$$
In the notation $\mu^*=\mathcal{L}(X,\widetilde{X})$ and $\mu_0^*=\mathcal{L}(T,\widetilde{T})$, it suffices to let
$$\mu_0^*(A)=\int_{\mathbb{R}^{2p}}\mathcal{N}_{2p}(x,\,c I)(A)\,\mu^*(dx)\quad\quad\text{for all Borel sets }A\subset\mathbb{R}^{2p}$$
where $c$ is a suitable constant. Thus, $\mathcal{L}(X,\widetilde{X})$ can be approximated in total variation by $\mathcal{L}(T,\widetilde{T})$ for any knockoff $\widetilde{X}$ which makes $\mathcal{L}(X,\widetilde{X})$ absolutely continuous with a Lipschitz density. While this fact is theoretically meaningful and supports the CIKs further, the above formula for $\mu_0^*$ has little practical use, since $\mu^*$ is generally unknown (it is even unknown how to obtain $\widetilde{X}$).

\subsection{A characterization of the CIKs}\label{e3hhu7}
Recall that
$$\mathcal{K}=\bigl\{\widetilde{X}:\,\mathcal{L}(X,\widetilde{X})\text{ admits representation \eqref{g4rx3e7} for some }\Theta,\,\gamma\text{ and }P_i(\cdot\mid\theta)\bigr\}$$
is the class of the CIKs of $X$. Such a $\mathcal{K}$ does not include all possible knockoffs. Here is a trivial example.

\begin{ex}\label{u8j3s}\textbf{(Not every knockoff is a CIK).}
Suppose that $X_1,\ldots,X_p$ are i.i.d. with $P(X_1=0)=P(X_1=1)=1/2$. In this case, it would be natural to take $\widetilde{X}$ as an independent copy of $X$. But suppose we let
$$\widetilde{X}=(\widetilde{X}_1,\ldots,\widetilde{X}_p)=(1-X_1,\ldots,1-X_p).$$
Then, for all $a,\,b\in\{0,1\}^p$,
$$P(X=a,\widetilde{X}=b)=P(X=a)\quad\text{if }b_i=1-a_i\text{ for each }i=1,\ldots,p$$
while $P(X=a,\widetilde{X}=b)=0$ otherwise. Based on this fact, it is straightforward to verify that $\widetilde{X}$ is a knockoff copy of $X$. However, since $X_i^2=X_i$, one obtains
$$\text{cov}(X_i,\widetilde{X}_i)=E\bigl\{X_i(1-X_i)\bigr\}-E(X_i)^2=E(X_i)-E(X_i^2)-E(X_i)^2=-E(X_i)^2<0.$$
Now, if $\widetilde{X}\in\mathcal{K}$, Jensen's inequality implies $\text{cov}(X_i,\widetilde{X}_i)\ge 0$. Hence, $\widetilde{X}\notin\mathcal{K}$.
\end{ex}

\medskip

\noindent Based on Example \ref{u8j3s}, a question is how to identify the members of $\mathcal{K}$ among all possible knockoffs $\widetilde{X}$. To answer this question, we recall that $(X,\widetilde{X})$ is said to be infinitely extendable if there exists an (infinite) sequence $V=(V_1,V_2,\ldots)$ such that $(V_1,\ldots,V_{2p})\sim (X,\widetilde{X})$ and $V$ satisfies the same invariance condition as $(X,\widetilde{X})$. Formally, the latter request should be meant as follows. Given three integers $i,\,j,\,k$ with $1\le i\le p$ and $j,\,k\ge 0$, define a new sequence $V^*=(V_1^*,V_2^*,\ldots)$ by swapping $V_{kp+i}$ with $V_{jp+i}$ and leaving all other elements of $V$ fixed, that is,
$$ V^*_{kp+i}=V_{jp+i},\quad V^*_{jp+i}=V_{kp+i},\quad V^*_r=V_r\,\text{ if }r\notin\{kp+i,\,jp+i\}.$$
Then, $V$ is required to satisfy
\begin{gather}\label{h5x2w}
V^*\sim V\quad\quad\text{for all }1\le i\le p\text{ and }j,\,k\ge 0.
\end{gather}
Condition \eqref{h5x2w} is nothing but a form of partial exchangeability; see \cite{ALDOUS} and \cite{DIAFRE}. In fact, the main tool for proving the next result is de Finetti's theorem for partially exchangeable sequences.

\begin{thm}\label{t3}
Let $\widetilde{X}$ be a knockoff copy of $X$. Then, $\widetilde{X}\in\mathcal{K}$ if and only if $(X,\widetilde{X})$ is infinitely extendable.
\end{thm}

\medskip

\noindent The essence of Theorem \ref{t3} is that, if $(X,\widetilde{X})$ is required to be infinitely extendable, then $X$ {\em must} be conditionally independent (given some $Z$) and $\widetilde{X}$ {\em must} be the CIK of $X$. One reason for requiring infinite extendability has been given in Section \ref{6b8m3d0k}. Essentially, infinite extendability of $(X,\widetilde{X})$ amounts to conditional independence between $X$ and $\widetilde{X}$, which in turn implies optimality of $\widetilde{X}$ under various criterions for increasing the power of the knockoff procedure; see Example \ref{v63e0qa2}.

\medskip

\begin{proof}[Proof of Theorem \ref{t3}]
Suppose $\widetilde{X}\in\mathcal{K}$, that is, $\mathcal{L}(X,\widetilde{X})$ admits representation \eqref{g4rx3e7} for some $\Theta$, $\gamma$ and $P_i(\cdot\mid\theta)$. For all $A_1,A_2,\ldots\in\mathcal{B}$, define
$$P(V_1\in A_1,V_2\in A_2,\ldots)=\int_\Theta\prod_{k=0}^\infty P_1(A_{kp+1}\mid\theta)\,P_2(A_{kp+2}\mid\theta)\ldots P_p(A_{kp+p}\mid\theta)\,\gamma(d\theta).$$
Then, $V=(V_1,V_2,\ldots)$ is an infinite sequence satisfying condition \eqref{h5x2w}. Moreover, by \eqref{g4rx3e7}, one obtains
\begin{gather*}
P(V_1\in A_1,\ldots,V_p\in A_p,V_{p+1}\in B_1,\ldots,V_{2p}\in B_p)
\\=\int_\Theta P_1(A_1\mid\theta)\ldots P_p(A_p\mid\theta)\,P_1(B_1\mid\theta)\ldots P_p(B_p\mid\theta)\,\gamma(d\theta)
\\=P(X_1\in A_1,\ldots,X_p\in A_p,\widetilde{X}_1\in B_1,\ldots,\widetilde{X}_p\in B_p)
\end{gather*}
whenever $A_i,\,B_i\in\mathcal{B}$ for each $i$. Hence, $(X,\widetilde{X})$ is infinitely extendable. Conversely, suppose $(X,\widetilde{X})$ is infinitely extendable and take an infinite sequence $V=(V_1,V_2,\ldots)$ satisfying condition \eqref{h5x2w} and $(V_1,\ldots,V_{2p})\sim (X,\widetilde{X})$. Let $\mathcal{Q}$ denote the set of all probability measures on $\mathbb{R}$. By \eqref{h5x2w}, $V$ is partially exchangeable; see e.g. \cite{ALDOUS}. Hence, by de Finetti's theorem, there is a probability measure $\pi$ on $\mathcal{Q}^p$ such that
\begin{gather*}
P(V_1\in A_1,V_2\in A_2,\ldots)=\int_{\mathcal{Q}^p}\prod_{k=0}^\infty q_1(A_{kp+1})\,q_2(A_{kp+2})\ldots q_p(A_{kp+p})\,\pi(dq_1,\ldots,dq_p)
\end{gather*}
for all $A_1,A_2,\ldots\in\mathcal{B}$; see \cite{ALDOUS} again. (Such a $\pi$ is usually called the de Finetti's measure of $V$). In particular,
\begin{gather*}
P(X_1\in A_1,\ldots,X_p\in A_p,\widetilde{X}_1\in B_1,\ldots,\widetilde{X}_p\in B_p)
\\=P(V_1\in A_1,\ldots,V_p\in A_p,V_{p+1}\in B_1,\ldots,V_{2p}\in B_p)
\\=\int_{\mathcal{Q}^p}q_1(A_1)\ldots q_p(A_p)\,q_1(B_1)\ldots q_p(B_p)\,\pi(dq_1,\ldots,dq_p).
\end{gather*}
Thus, to conclude the proof, it suffices to let
$$\Theta=\mathcal{Q}^p,\quad\gamma=\pi,\quad P_i(\cdot\mid\theta)=q_i(\cdot)\text{ for all }\theta=(q_1,\ldots,q_p)\in\mathcal{Q}^p.$$
\end{proof}

\section{2-valued and 3-valued covariates}\label{z33w21h}

In applications, an important special case is $X_i\in\{0,1\}$. In a genetic framework, for instance, $X_i=0$ or $X_i=1$ according to whether the $i$-th gene is absent or present. Another meaningful case is $X_i\in\{0,1,2\}$, where $X_i=2$ can be given various interpretations. For instance, $X_i=2$ could mean that the absence/presence of the $i$-th gene cannot be established. Despite their practical significance, to our knowledge, these cases have not received much attention, from the theoretical point of view, to date. In this section, we try to fill this gap. We aim to build a CIK $\widetilde{X}$ when $X$ is a vector of 2-valued or 3-valued random variables.

\medskip

\noindent There are obviously various cases. For instance, some covariates are 2-valued, other 3-valued, and the remaining ones have a continuous distribution function. Here, we only focus on two extreme situations: either all covariates are 2-valued or all are 3-valued.

\medskip

\noindent Suppose first $X_i\in\{0,1\}$ for all $i$. To build a CIK, $X_1,\ldots,X_p$ must be conditionally independent given a random parameter $\theta$. Here, it is natural to let $\theta=(\theta_1,\ldots,\theta_p)$ with $\theta_i\in (0,1)$ regarded as the (random) probability of the event $\{X_i=1\}$. Accordingly, $X_1,\ldots,X_p$ are assumed conditionally independent given $\theta=(\theta_1,\ldots,\theta_p)$ with $P(X_i=1\mid\theta)=\theta_i$. In this case,
$$P(X=x)=P(X_1=x_1,\ldots,X_p=x_p)=\int_{(0,1)^p}\prod_{i=1}^p\theta_i^{x_i}(1-\theta_i)^{1-x_i}\,\gamma(d\theta)$$
where $\gamma$ denotes the probability distribution of $\theta$ and
$$x=(x_1,\ldots,x_p)\in\{0,1\}^p.$$
To be concrete, we also assume
$$\theta=\lambda\,c$$
where $\lambda\in (0,1)$ is a random scalar and $c=(c_1,\ldots,c_p)\in (0,1)^p$ a vector of known constants. Moreover, we take $\lambda$ uniformly distributed on $(0,1)$ and we let
$$s=\sum_{i=1}^px_i,\quad S=\bigl\{i:x_i=1\bigr\},\quad d_0=1,\quad d_j=\sum_{\begin{array}{c}
                                                                             i_1<\ldots<i_j \\
                                                                                 i_1,\ldots,i_j\notin S \\
                                                                               \end{array}} c_{i_1}\,c_{i_2}\ldots c_{i_j}.$$
Then, after some algebra, one obtains
\begin{gather*}
P(X=x)=\prod_{i\in S}c_i\,\int_0^1\lambda^s\,\prod_{i\notin S}(1-\lambda\,c_i)\,d\lambda=\prod_{i\in S}c_i\,\,\sum_{j=0}^{p-s}(-1)^j\,\frac{d_j}{j+s+1}.
\end{gather*}
Similarly, we can evaluate $P(X=x,\,\widetilde{X}=\widetilde{x})$ where
$$\widetilde{x}=(\widetilde{x}_1,\ldots,\widetilde{x}_p)\in\{0,1\}^p.$$
To this end, define
$$t=\sum_{i=1}^p\widetilde{x}_i,\quad T=\bigl\{i:\widetilde{x}_i=1\bigr\},\quad e_0=1,\quad e_j=\sum_{\begin{array}{c}
                                                                             i_1<\ldots<i_j \\
                                                                                 i_1,\ldots,i_j\notin T \\
                                                                               \end{array}} c_{i_1}\,c_{i_2}\ldots c_{i_j}.$$
Then,
\begin{gather*}
P(X=x,\,\widetilde{X}=\widetilde{x})=\prod_{i\in S}c_i\,\prod_{i\in T}c_i\,\int_0^1\lambda^{s+t}\,\prod_{i\notin S}(1-\lambda\,c_i)\,\prod_{i\notin T}(1-\lambda\,c_i)\,d\lambda
\\=\prod_{i\in S}c_i\,\prod_{i\in T}c_i\,\,\sum_{j=0}^{p-s}\sum_{k=0}^{p-t}(-1)^{j+k}\,\frac{d_je_k}{j+k+s+t+1}.
\end{gather*}
Finally,
\begin{gather*}
P\bigl(\widetilde{X}=\widetilde{x}\mid X=x\bigr)=\frac{P\bigl(X=x,\,\widetilde{X}=\widetilde{x}\bigr)}{P(X=x)}
\\=\prod_{i\in T}c_i\,\,\left(\sum_{j=0}^{p-s}(-1)^j\,\frac{d_j}{j+s+1}\right)^{-1}\,\,\sum_{j=0}^{p-s}\sum_{k=0}^{p-t}(-1)^{j+k}\,\frac{d_je_k}{j+k+s+t+1}.
\end{gather*}

\medskip

\noindent We now have an explicit formula for $\mathcal{L}(\widetilde{X}\mid X=x)$. In a sense, this is the best we can do. In fact, after observing $X=x$, a value $\widetilde{x}$ for $\widetilde{X}$ can be drawn directly from $\mathcal{L}(\widetilde{X}\mid X=x)$.

\medskip

\noindent Next, suppose that $X_i\in\{0,1,2\}$ for all $i$. To deal with this case, we assume $X_1,\ldots,X_p$ conditionally independent given $\lambda$ with
$$P(X_i=0\mid\lambda)=\lambda\,(1-c_i),\quad P(X_i=1\mid\lambda)=\lambda\,c_i,\quad P(X_i=2\mid\lambda)=1-\lambda,$$
where $\lambda\in (0,1)$ is a random scalar and $c_i\in (0,1)$ a fixed known constant. We give $\lambda$ a beta distribution with parameters $a>0$ and $b>0$. Moreover, for all
$$x=(x_1,\ldots,x_p)\in\{0,1,2\}^p\quad\text{and}\quad\widetilde{x}=(\widetilde{x}_1,\ldots,\widetilde{x}_p)\in\{0,1,2\}^p,$$
we let
\begin{gather*}
S_0=\bigl\{i:x_i=0\bigr\},\quad S_1=\bigl\{i:x_i=1\bigr\},\quad T_0=\bigl\{i:\widetilde{x}_i=0\bigr\},\quad T_1=\bigl\{i:\widetilde{x}_i=1\bigr\},
\\m_2=\text{card}\,\bigl\{i:x_i=2\bigr\}\quad\text{and}\quad n_2=\text{card}\,\bigl\{i:\widetilde{x}_i=2\bigr\}.
\end{gather*}
Then,
\begin{gather*}
P(X=x)=\int_0^1\prod_iP(X_i=x_i\mid\lambda)\,\,\frac{\Gamma(a+b)}{\Gamma(a)\Gamma(b)}\,\lambda^{a-1}(1-\lambda)^{b-1}d\lambda
\\=\frac{\Gamma(a+b)}{\Gamma(a)\Gamma(b)}\,\,\prod_{i\in S_1}c_i\,\prod_{i\in S_0}(1-c_i)\,\,\int_0^1\lambda^{a+p-m_2-1}(1-\lambda)^{b+m_2-1}d\lambda
\\=\frac{\Gamma(a+b)}{\Gamma(a)\Gamma(b)}\,\,\frac{\Gamma(a+p-m_2)\,\Gamma(b+m_2)}{\Gamma(a+b+p)}\,\prod_{i\in S_1}c_i\,\prod_{i\in S_0}(1-c_i)
\end{gather*}
and
\begin{gather*}
P\bigl(X=x,\,\widetilde{X}=\widetilde{x}\bigr)=\int_0^1\prod_iP(X_i=x_i\mid\lambda)\,P(X_i=\widetilde{x}_i\mid\lambda)\,\,\frac{\Gamma(a+b)}{\Gamma(a)\Gamma(b)}\,\lambda^{a-1}(1-\lambda)^{b-1}d\lambda
\\=\frac{\Gamma(a+b)}{\Gamma(a)\Gamma(b)}\,\,\frac{\Gamma(a+2p-m_2-n_2)\,\Gamma(b+m_2+n_2)}{\Gamma(a+b+2p)}\,\prod_{i\in S_1}c_i\,\prod_{i\in T_1}c_i\,\prod_{i\in S_0}(1-c_i)\,\prod_{i\in T_0}(1-c_i).
\end{gather*}
Hence, even in this case, we have an explicit formula for $\mathcal{L}(\widetilde{X}\mid X=x)$, that is
\begin{gather*}
P\bigl(\widetilde{X}=\widetilde{x}\mid X=x\bigr)=\frac{P\bigl(X=x,\,\widetilde{X}=\widetilde{x}\bigr)}{P(X=x)}
\\=\frac{\Gamma(a+2p-m_2-n_2)\,\Gamma(b+m_2+n_2)\,\Gamma(a+b+p)}{\Gamma(a+b+2p)\,\Gamma(a+p-m_2)\,\Gamma(b+m_2)}\,\,\prod_{i\in T_1}c_i\,\prod_{i\in T_0}(1-c_i).
\end{gather*}

\section{Mixtures of centered normal distributions}\label{mcnd555}

In this section, $\theta=(\theta_1,\ldots,\theta_p)$ is a vector of strictly positive random variables and $X_1,\ldots,X_p$ are conditionally independent given $\theta$ with
$$X_i\mid\theta\,\sim\,\mathcal{N}_1(0,\theta_i)\quad\quad\text{for each }i=1,\ldots,p.$$

\medskip

\noindent Mixtures of centered normal distributions allow to model various real situations while preserving some properties of the Gaussian laws. For this reason, they are quite popular in applications; see e.g. \cite{HHL} and references therein. Among other things, they are involved in Bayesian inference for logistic models \cite{PSW} and they arise as the limit laws in the CLT for exchangeable random variables \cite[Sect. 3]{BPR04}. A further motivation for this type of data is that $E(X_i\mid\theta)=0$. Hence, the CIKs are optimal and in particular
$$\text{cov}(X_i,\widetilde{X}_i)=0\quad\quad\text{for all }i;$$
see Example \ref{v63e0qa2}.

\medskip

\noindent To build a CIK, a "prior" $\gamma$ on $\Theta=(0,\infty)^p$ is to be selected. Quite surprisingly, to our knowledge, the choice of $\gamma$ seems to be almost neglected in the Bayesian
literature (apart from the special case $\theta_1=\ldots=\theta_p$); see e.g. \cite{GCSDVR}. We next propose two choices of $\gamma$. As in Section \ref{z33w21h}, we let $\theta=\lambda\,c$ where $\lambda>0$ is a scalar and $c=(c_1,\ldots,c_p)$ a vector such that $c_i>0$ for all $i$.

\subsection{First choice of $\gamma$}\label{h78uj} We first assume that $\lambda$ is random but $c$ is not. Equivalently, we suppose that the ratios $\theta_i/\theta_j=c_i/c_j$ are non-random and known. While simple, this assumption makes sense in various applications, for instance in a financial framework.

\medskip

\noindent The random variable $\lambda$ is given an inverse Gamma distribution with parameters $a>0$ and $b>0$, that is, $\lambda$ has density $\psi(x)=b^a\,\Gamma(a)^{-1}x^{-a-1}\exp{(-b/x)}$ for $x>0$. In this case, the density of $(X,\widetilde{X})$ is
\begin{gather*}
f(x,\widetilde{x})=\int_0^\infty\prod_{i=1}^p\,f_i(x_i\mid\lambda,c)\,f_i(\widetilde{x}_i\mid\lambda,c)\,\psi(\lambda)\,d\lambda
\end{gather*}
where $x$ and $\widetilde{x}$ are points of $\mathbb{R}^p$ and $f_i(\cdot\mid\lambda,c)$ is the density of $\mathcal{N}_1(0,\,c_i\,\lambda)$. Hence,
\begin{gather*}
f(x,\widetilde{x})=\frac{b^a}{(2\pi)^p\,\Gamma(a)\,\prod_{i=1}^pc_i}\,\,\int_0^\infty\lambda^{-a-p-1}\exp\Bigl\{-\frac{1}{\lambda}\,\Bigl(b+\sum_{i=1}^p\frac{x_i^2+\widetilde{x}_i^2}{2c_i}\Bigr)\Bigr\}\,d\lambda
\\=\frac{b^a\,\Gamma(a+p)}{(2\pi)^p\,\Gamma(a)\,\prod_{i=1}^pc_i}\,\,\Bigl(b+\sum_{i=1}^p\frac{x_i^2+\widetilde{x}_i^2}{2c_i}\Bigr)^{-(a+p)}.
\end{gather*}
Similarly, the density of $X$ is
\begin{gather*}
h(x)=\int_0^\infty\prod_{i=1}^p\,f_i(x_i\mid\lambda,c)\,\psi(\lambda)\,d\lambda=\frac{b^a\,\Gamma(a+p/2)}{(2\pi)^{p/2}\,\Gamma(a)\,\sqrt{\prod_{i=1}^pc_i}}\,\,\Bigl(b+\sum_{i=1}^p\frac{x_i^2}{2c_i}\Bigr)^{-(a+p/2)}.
\end{gather*}

\medskip

\noindent It is worth noting that $f$ and $h$ are densities of Student's-$t$ distributions. We recall that the $m$-variate Student's-$t$ distribution with $k$ degrees of freedom is the absolutely continuous distribution on $\mathbb{R}^m$ with density
$$\varphi(x)=\frac{\Gamma\bigl[(m+k)/2\bigr]}{\Gamma(k/2)\,(k\pi)^{m/2}\,\sqrt{\text{det}\Sigma}}\,\,\Bigl(1+(1/k)\,x^T\Sigma^{-1}x\Bigr)^{-(m+k)/2}\quad\text{for each }x\in\mathbb{R}^m,$$
where $\Sigma$ is a symmetric positive definite $m\times m$ matrix. Hence, one obtains $\varphi=f$ if $m=2p$, $k=2a$ and $\Sigma=b\,a^{-1}\,$diag$(c_1,\ldots,c_p,c_1,\ldots,c_p)$ and $\varphi=h$ if $m=p$, $k=2a$ and $\Sigma=b\,a^{-1}\,$diag$(c_1,\ldots,c_p)$.

\medskip

\noindent Finally, the conditional density of $\widetilde{X}$ given $X=x$ can be written as
$$g(\widetilde{x}\mid x)=\frac{f(x,\widetilde{x})}{h(x)}=\frac{\Gamma(a+p)}{(2\pi)^{p/2}\,\Gamma(a+p/2)\,\sqrt{\prod_{i=1}^pc_i}}\,\,\,\frac{\Bigl(b+\sum_{i=1}^p\frac{x_i^2}{2c_i}\Bigr)^{a+p/2}}{\Bigl(b+\sum_{i=1}^p\frac{x_i^2+\widetilde{x}_i^2}{2c_i}\Bigr)^{a+p}}.$$
Once again, $g(\cdot\mid x)$ is the density of a Student's-$t$ distribution (with parameters depending on $x$). To see this, it suffices to let $m=p$, $k=2a+p$, and
$$\Sigma=\frac{2}{2a+p}\,\Bigl(b+\sum_{i=1}^p\frac{x_i^2}{2c_i}\Bigr)\,\text{diag}(c_1,\ldots,c_p).$$
Thus, we have an explicit formula for $g(\cdot\mid x)$ and this is quite useful in applications. A numerical example is in Section \ref{n5e}.

\subsection{Second choice of $\gamma$}\label{a2l0h5m} Suppose now that $c$ is random and independent of $\lambda$. Let $c$ be given an absolutely continuous distribution with density $q$. Then, $f$, $h$ and $g$ turn into
\begin{gather*}
f(x,\widetilde{x})=\int_{(0,\infty)^p}\int_0^{+\infty}\prod_{i=1}^p\,f_i(x_i\mid\lambda,c)\,f_i(\widetilde{x}_i\mid\lambda,c)\,\psi(\lambda)\,q(c)\,d\lambda \, dc,
\\h(x)=\int_{(0,\infty)^p}\int_0^{+\infty}\prod_{i=1}^p\,f_i(x_i\mid\lambda,c)\,\psi(\lambda)\,q(c)\,d\lambda \, dc\quad\text{and}\quad g(\widetilde{x}\mid x)=\frac{f(x,\widetilde{x})}{h(x)}.
\end{gather*}
As an example, $c_1,\ldots,c_p$ could be taken i.i.d. according to a uniform distribution on some bounded interval $B\subset (0,\infty)$, i.e.,
$$q(c)=\frac{1}{\text{length}(B)^p}\,\prod_{i=1}^p1_B(c_i).$$

\medskip

\noindent In general, the above integrals cannot be explicitly evaluated. Hence, sampling from $g(\cdot\mid x)$ is not easy, but it is still possible by computational methods. For instance, we could proceed as follows. Since $g(\cdot\mid x)$ is proportional to $f(x,\cdot)$, we focus on $f(x,\cdot)$. Then, to sample from $f(x,\cdot)$, we adopt a data augmentation strategy where $\lambda$ and $c$ are treated as auxiliary variables. The idea is to consider the density function
$$g^*(\widetilde{x},\lambda,c\mid x)\propto \prod_{i=1}^p\,f_i(x_i\mid\lambda,c)\,f_i(\widetilde{x}_i\mid\lambda,c)\,\psi(\lambda)\,q(c)$$
and perform a Gibbs sampling on the variables $(\widetilde{x},\lambda,c)$.  We conclude this section by listing the full conditional distributions required to run the algorithm.

\medskip

\begin{itemize}
\item Let $\widetilde{x}_{-i}=(\widetilde{x}_1,\ldots,\widetilde{x}_{i-1},\widetilde{x}_{i+1},\ldots,\widetilde{x}_p)$. The full conditional distribution of $(\widetilde{x}_i\mid \widetilde{x}_{-i},\lambda,c)$ is proportional to $f_i(\cdot\mid\lambda,c)$. This means that
$\widetilde{x}_i$ can be sampled from a centered normal distribution with variance $\lambda c_i$.

\medskip

\item The full conditional distribution of $(\lambda\mid \widetilde{x}, c)$ is proportional to

$$\frac{\psi(\lambda)}{\lambda^p}\,\exp\Bigl\{-\frac{1}{\lambda}\,\sum_{i=1}^p\frac{x_i^2+\widetilde{x}_i^2}{2c_i}\Bigr\}.$$

Hence, since $\psi$ is the inverse gamma density with parameters $a$ and $b$, the full conditional of $\lambda$ is still an inverse gamma with parameters
$$a^*=a+p\quad\text{and}\quad b^*=b+\sum_{i=1}^p\frac{x_i^2+\widetilde{x}_i^2}{2c_i}.$$
Obviously, $\lambda$ could be also given a different distribution. In this case, the corresponding full conditional is probably more involved, but one may use a metropolis within Gibbs step.

\medskip

\item Let $c_{-i}=(c_1,\ldots,c_{i-1},c_{i+1},\ldots,c_p)$. The full conditional distribution of $(c_i\mid \widetilde{x},\lambda,c_{-i})$ is proportional to

$$\frac{q(c)}{c_i}\,\exp\Bigl\{-\frac{1}{c_i}\,\frac{x_i^2+\widetilde{x}_i^2}{2\lambda}\Bigr\}.$$

Sampling from the above is not straightforward and may require a metropolis within Gibbs step.

\end{itemize}

\section{A numerical experiment}\label{n5e}

In this section, the CIKs are tested numerically against both simulated and real data. To this end, $X$ is assumed to be as in Section \ref{h78uj}. Hence, $\mathcal{L}(X)$ is a mixture of centered normal distributions and $\theta=\lambda\,c$, where the scalar $\lambda$ has an inverse gamma distribution with parameters $a$ and $b$ while $c=(c_1,\ldots,c_p)$ is a vector of strictly positive constants.

\medskip

\noindent To learn something about the impact of the parameters, the experiment has been repeated for various choices of $a$, $b$ and $c$. The obtained results are quite stable with respect to $a$ and $b$ but exhibit a notable variability with respect to $c$. In the sequel, $a$ and $b$ have been selected so as to control the mean and the variance of $\lambda$ (which hold $b/(a-1)$ and $\frac{b^2}{(a-1)^2(a-2)}$, respectively, for $a>2$). In case of real data (Section \ref{jc387hu}) $a$ and $b$ have been also tuned based on the observed value of $X$. The choice of $c$ is certainly more delicate. As in Section \ref{a2l0h5m}, one option could be modeling $c$ as a random vector (rather than a fixed vector). For instance, $c_1,\ldots,c_p$ could be i.i.d, according to a uniform distribution on some interval, and independent of $\lambda$. However, in this section, $c$ is taken to be non-random. This choice has essentially three motivations. First, it may be convenient in real problems, in order to account for the different roles of the various covariates. Second, it is practically simpler since computational methods are not required. Third, if $c$ is non-random, a direct comparison with Section 6.3 of \cite{RSC} is easier.

\medskip

\noindent One more remark is in order. To compare different knockoff procedures, three popular criterions are the power, the false discovery rate, and the observed correlations between the $X_i$ and their knockoffs $\widetilde{X}_i$. However, as regards the CIKs of Section \ref{h78uj}, the third criterion is superfluos, since cov$(X_i,\widetilde{X}_i)=0$ for all $i$. Indeed, under the third criterion (as well as under the reconstructability criterion), the CIKs of Section \ref{h78uj} are superior to any other knockoff procedure; see Example \ref{v63e0qa2}.

\subsection{Simulated data}\label{s34rth8}

According to the usual format (see e.g. \cite{CFJL18} and \cite{RSC}) the simulation experiment has been performed as follows.

\begin{itemize}

\item A subset $I\subset\bigl\{1,\ldots,p\bigr\}$ such that $|I|=60$ has been randomly selected and the coefficients $\beta_1,\ldots,\beta_p$ have been defined as
$$\beta_i=0\,\text{ if }i\notin I\quad\text{and}\quad\beta_i=\frac{u}{\sqrt{n}}\,\text{ if }i\in I.$$
Here, $n$ is a positive integer and $u>0$ a parameter called {\em signal amplitude}.

\item $n$ i.i.d. observations
$$X^{(j)}=\bigl(X_{1j},\ldots,X_{pj}\bigr),\quad\quad j=1,\ldots,n,$$
have been generated from a $p$-variate Student's-$t$ distribution with $2a$ degrees of freedom and matrix $\Sigma=b\,a^{-1}\,$diag$(c_1,\ldots,c_p)$. Given $X^{(j)}$, the corresponding response variable $Y^{(j)}$ has been defined as
$$Y^{(j)}=\sum_{i=1}^p\beta_iX_{ij}+e_j$$
where $e_1,\ldots,e_n$ are i.i.d. standard normal errors.

\item For each $j=1,\ldots,n$, we sampled $m$ CIKs, say $\widetilde{X}^{(1,j)},\ldots,\widetilde{X}^{(m,j)}$, from the conditional distribution of $\widetilde{X}^{(j)}$ given $X^{(j)}=x^{(j)}$ where $x^{(j)}$ is the observed value of $X^{(j)}$. Precisely, for each $k=1,\ldots,m$, the value of $\widetilde{X}^{(k,j)}$ was sampled from the $p$-variate Student's-$t$ distribution with $2a+p$ degrees of freedom and matrix
$$\Sigma=\frac{2}{2a+p}\,\Bigl(b+\sum_{i=1}^p\frac{x_{ij}^2}{2c_i}\Bigr)\,\text{diag}(c_1,\ldots,c_p).$$

\item For each $k=1,\ldots,m$, the knockoff selection procedure has been applied to the data
$$\bigl\{Y^{(j)},\,X^{(j)},\,\widetilde{X}^{(k,j)}:j=1,\ldots,n\bigr\}$$
so as to calculate the power and the false discovery rate, say $pow(k)$ and $fdr(k)$. To do this, we exploited the R-cran package \verb"knockoff":

\verb"https://cran.r-project.org/web/packages/knockoff/index.html".

\noindent This package is based on the comparison between the lasso coefficient estimates of each covariate and its knockoff.

\item The final outputs are the arithmetic means of the powers and the false discovery rates, i.e.,
$$pow=(1/m)\,\sum_{k=1}^m pow(k)\quad\text{and}\quad fdr=(1/m)\,\sum_{k=1}^m fdr(k).$$

\end{itemize}

\medskip

\noindent To run the simulation experiment, we took $m=n=p=1000$ and a theoretical value of the false discovery rate equal to $0.1$. As already noted, the experiment has been repeated for various choices of the parameters $a,\,b,\,c,\,u$. Overall, the results have been quite stable with respect to all parameters but $c$.  The specific results reported here correspond to $a=6$, $b=10$, $c_i=i$ and $u=0.15, 0.2, 0.25, 0.3, 0.35, 0.4, 0.45, 0.5, 1, 1.5, 2, 2.5, 3$.

\medskip

\begin{figure}[!htbp]
        \centering
        \includegraphics[width=1\textwidth]{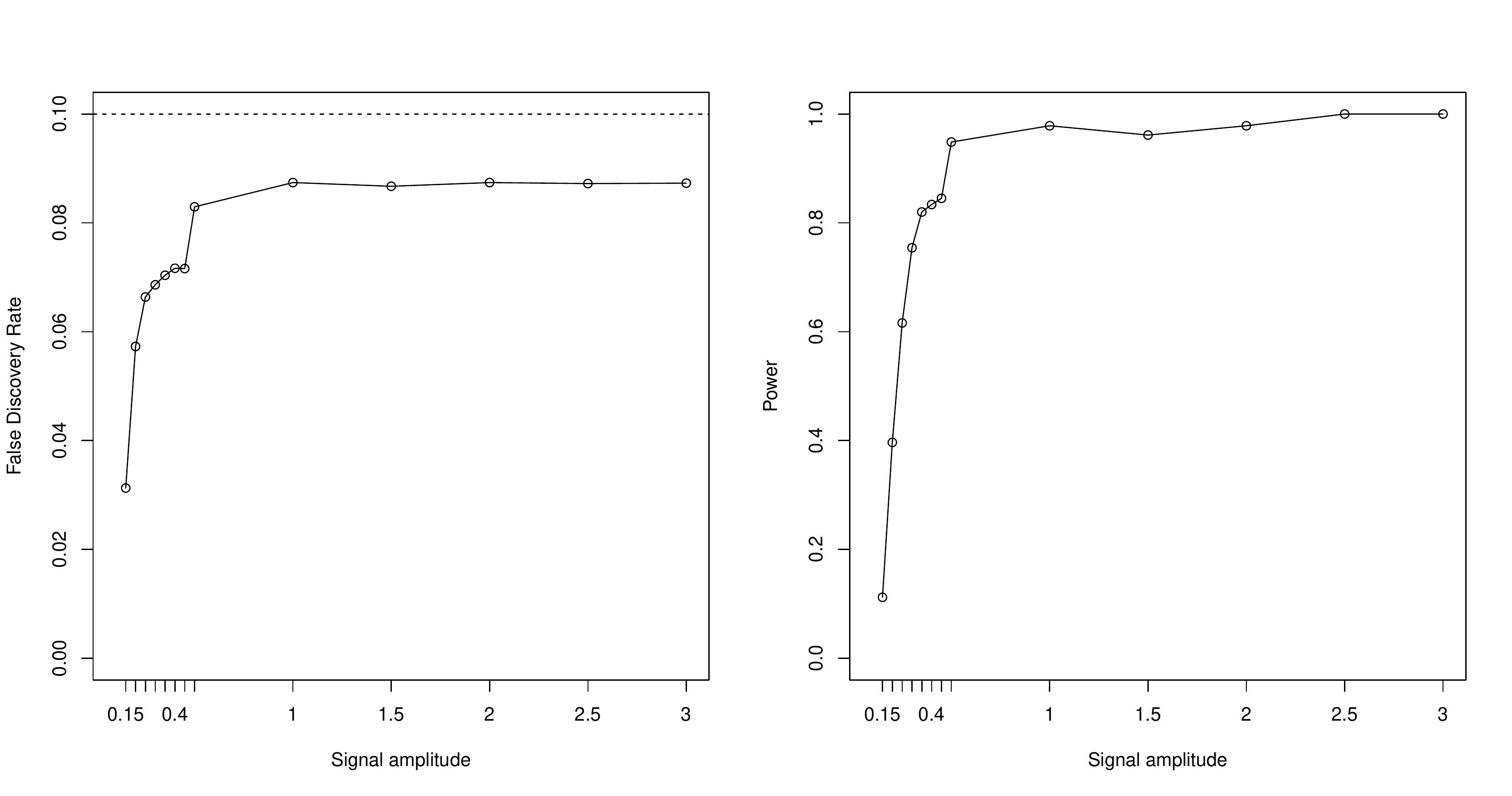}
        \caption{Results from the simulation experiment: False Discovery Rate (left) and Power (right)
performances for the CIKs with Signal amplitude equal to 0.15, 0.2, 0.25, 0.3, 0.35, 0.4, 0.45, 0.5, 1, 1.5, 2, 2.5, 3}\label{fig:simulation}
   \end{figure}

\noindent The observed results, in terms of $pow$ and $fdr$, are summarized in Figure \ref{fig:simulation}. The performance of the CIKs appears to be excellent, even if it slightly gets worse for small values of the amplitude $u$. It is worth noting that, as regards the power, the behavior of the CIKs is even optimal. This was quite expected, however, because of the optimality of the CIKs discussed in Example \ref{v63e0qa2}.

\medskip

\subsection{Real data}\label{jc387hu}

We next turn to real data. In this case, the CIKs can be compared with some other knockoff procedures, namely:  The Benjamin and Hochberg method \cite{BH95}, denoted by BHq; The fixed $X$ knockoff \cite{BC15}, denoted by Fixed-$X$; The model-$X$ Gaussian knockoff \cite{CFJL18}, denoted by Model-$X$; The second-order knockoff \cite{CFJL18,RSC}, denoted by Second-order. The comparison is based on the power, the false discovery rate, and the number of false and true discoveries. The results reported here correspond to $a=4$, $b=3$ and $c_i=i$.

\medskip

\noindent We focus on the human immunodeficiency virus type 1 (HIV-1) dataset \cite{HIV}, which has been used in several papers on the knockoff procedure; see e.g. \cite{BC15, RSC}. The dimension of our dataset is $n=846$ and $p=341$, where $n$ denotes the number of observations. The knockoff filter is applied to detect the mutations associated with drug resistance. In fact, the HIV-1 dataset provides drug resistance measurements. Furthermore, it includes genotype information from samples of HIV-1, with separate data sets for resistance to protease inhibitors, nucleoside reverses transcriptase inhibitors, and non-nucleoside RT inhibitors. We deal with resistance to protease inhibitors, and we analyze separately the following drugs: amprenavir (APV), atazanavir (ATV), indinavir (IDV), lopinavir (LPV), nelfinavir (NFV), ritonavir (RTV) and saquinavir (SQV).

\medskip

\begin{figure}[!htbp]
        \centering
        \includegraphics[width=1\textwidth]{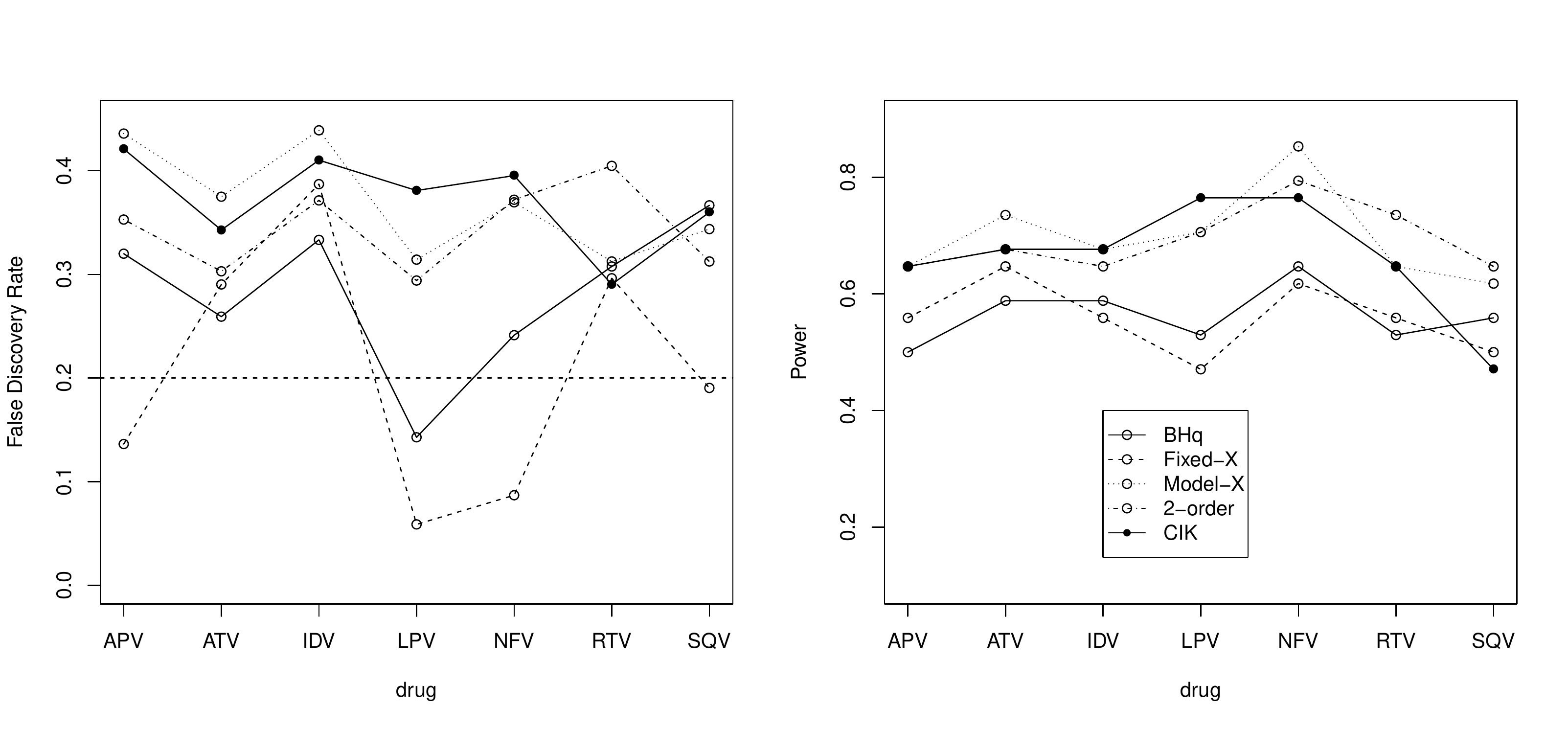}
          \caption{Results from the real example: False Discovery Rate (left) and Power (right) performances across methods and drugs} \label{fig:real}
\end{figure}

\begin{figure}[!htbp]
     \centering
     \hbox{
         \includegraphics[width=5cm]{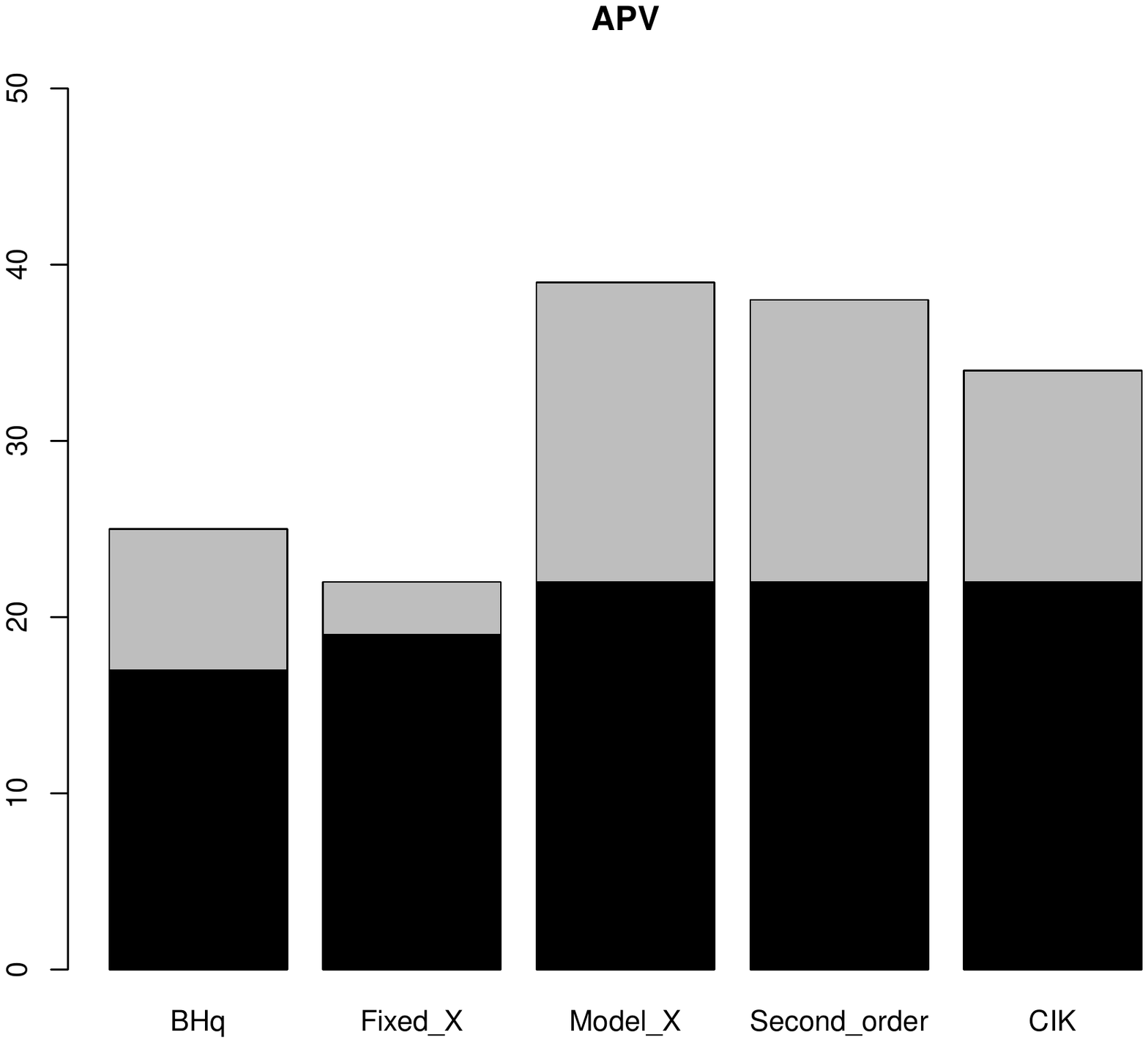}
         \includegraphics[width=5cm]{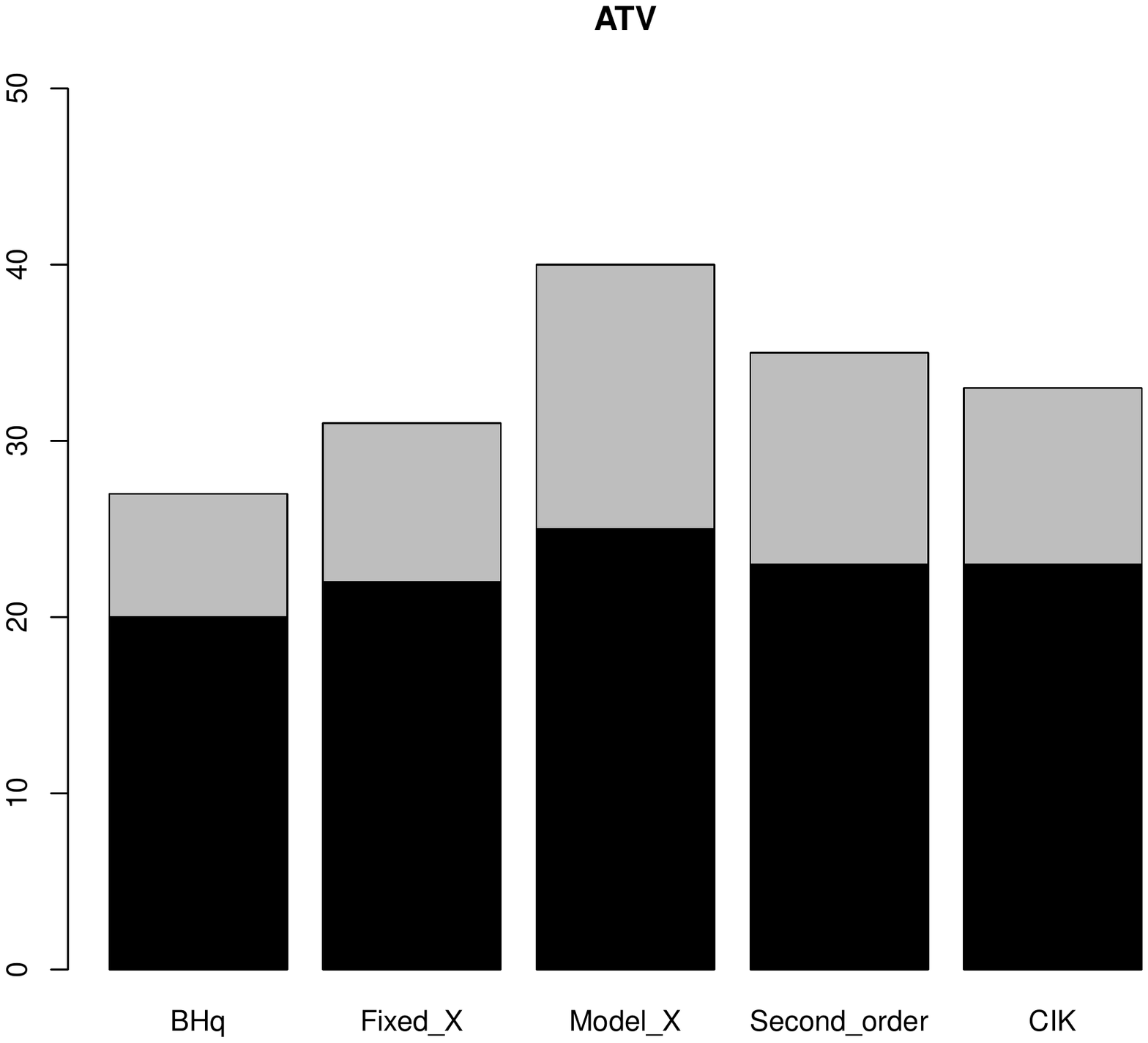}}
         \vspace{-3cm}
        \hbox{
         \includegraphics[width=5cm]{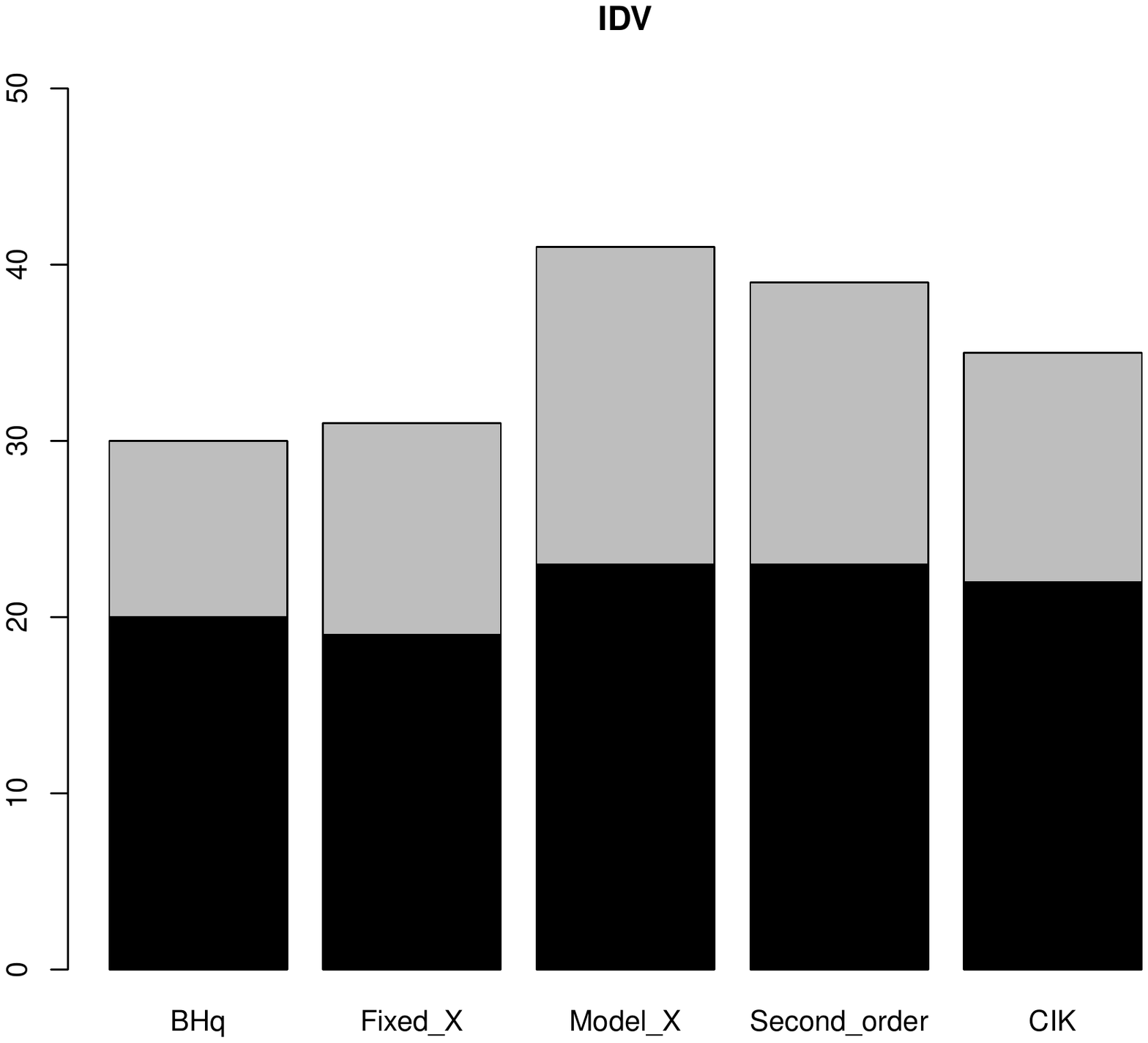}
         \includegraphics[width=5cm]{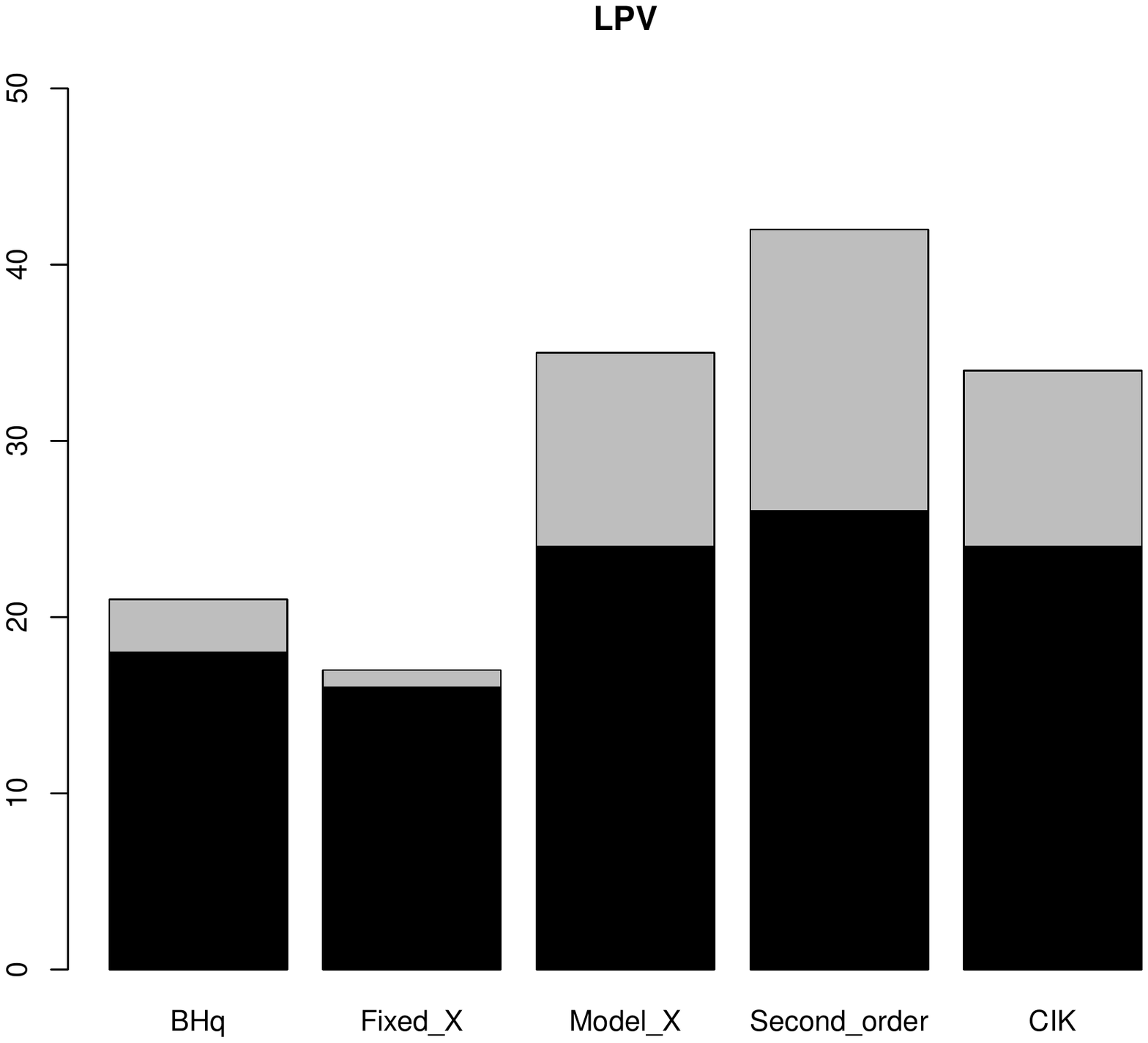}}
          \vspace{-3cm}
     \hbox{
         \includegraphics[width=5cm]{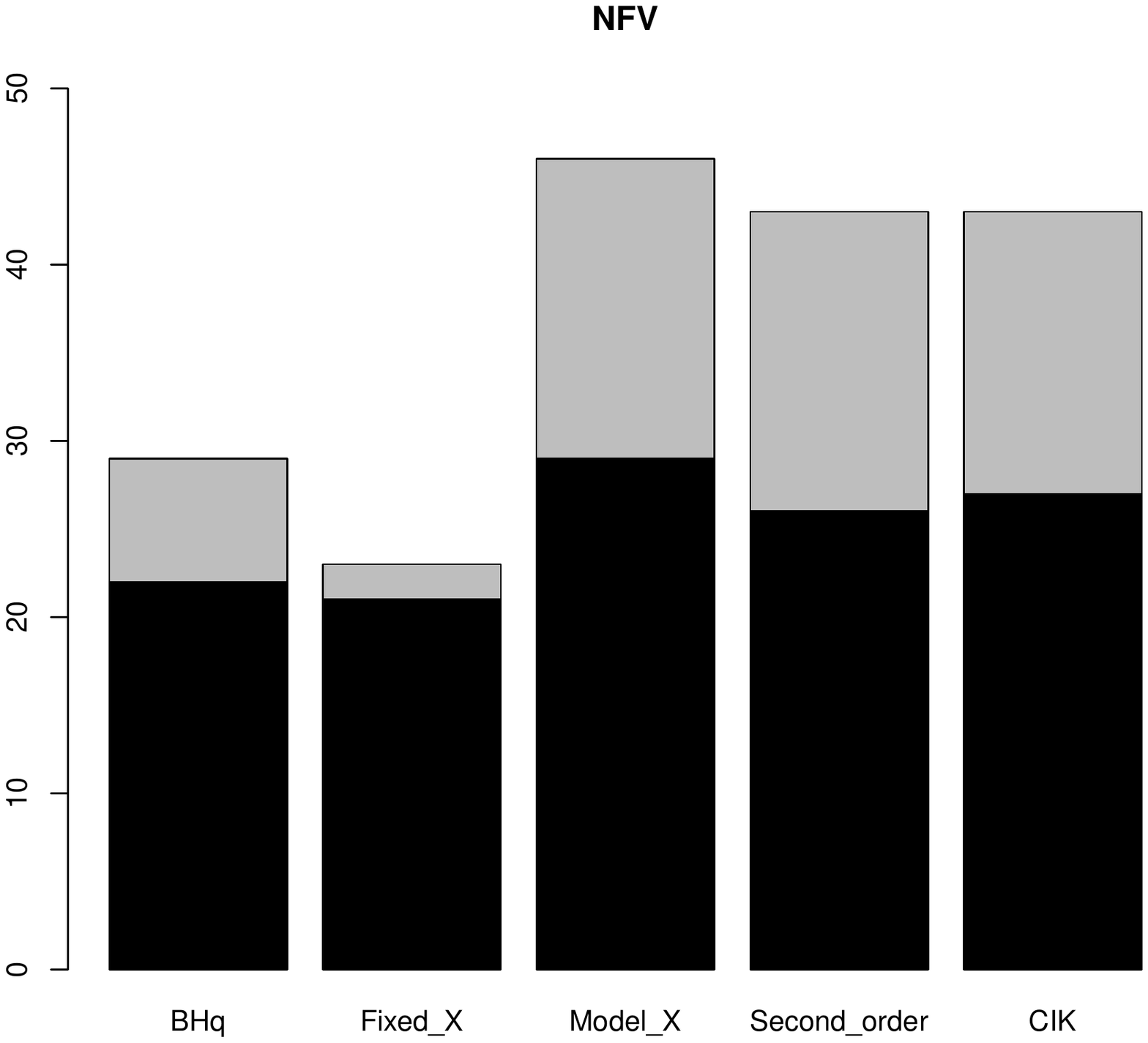}
         \includegraphics[width=5cm]{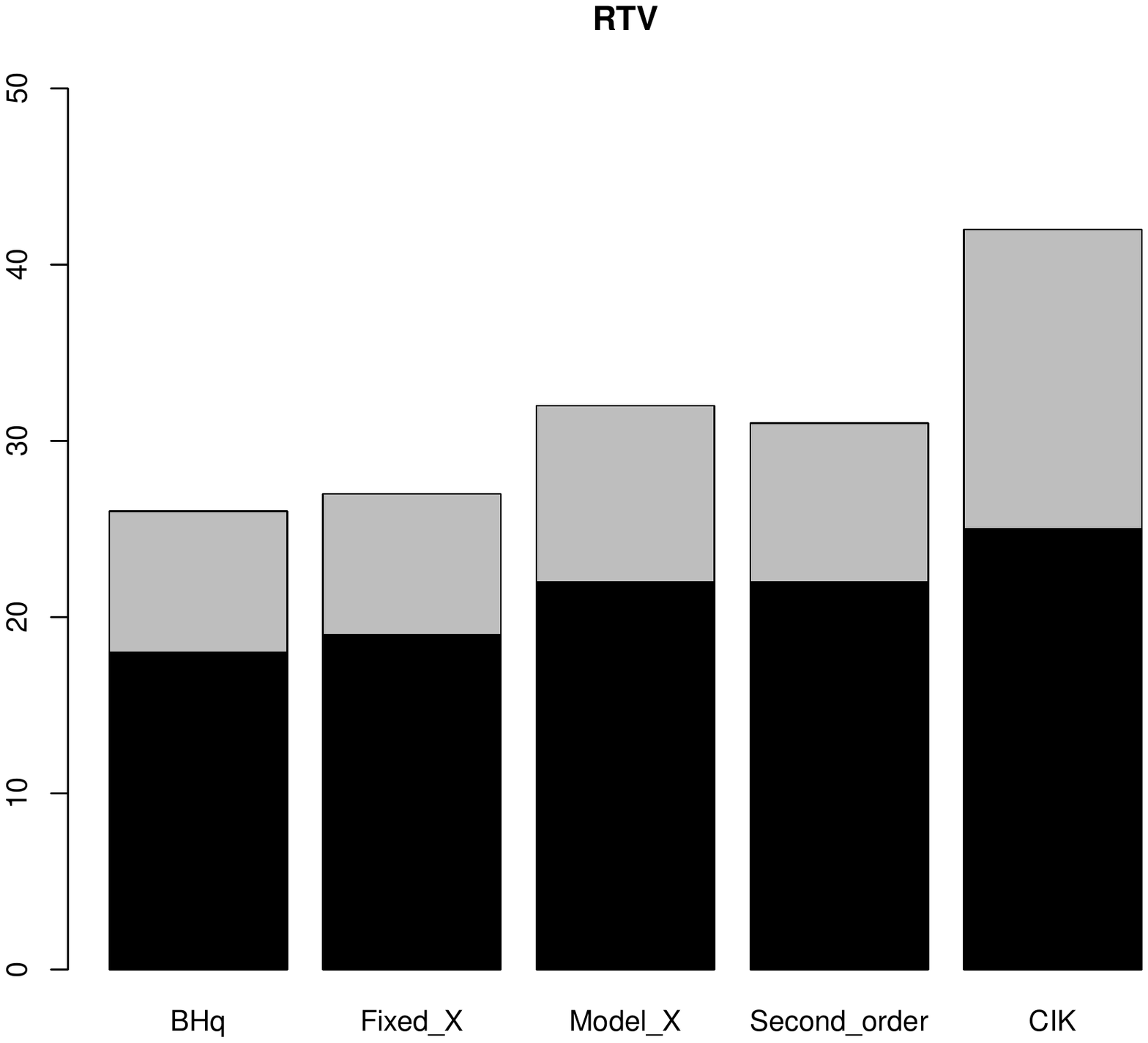}}
          \vspace{-3cm}
     \centering
         \includegraphics[width=5cm]{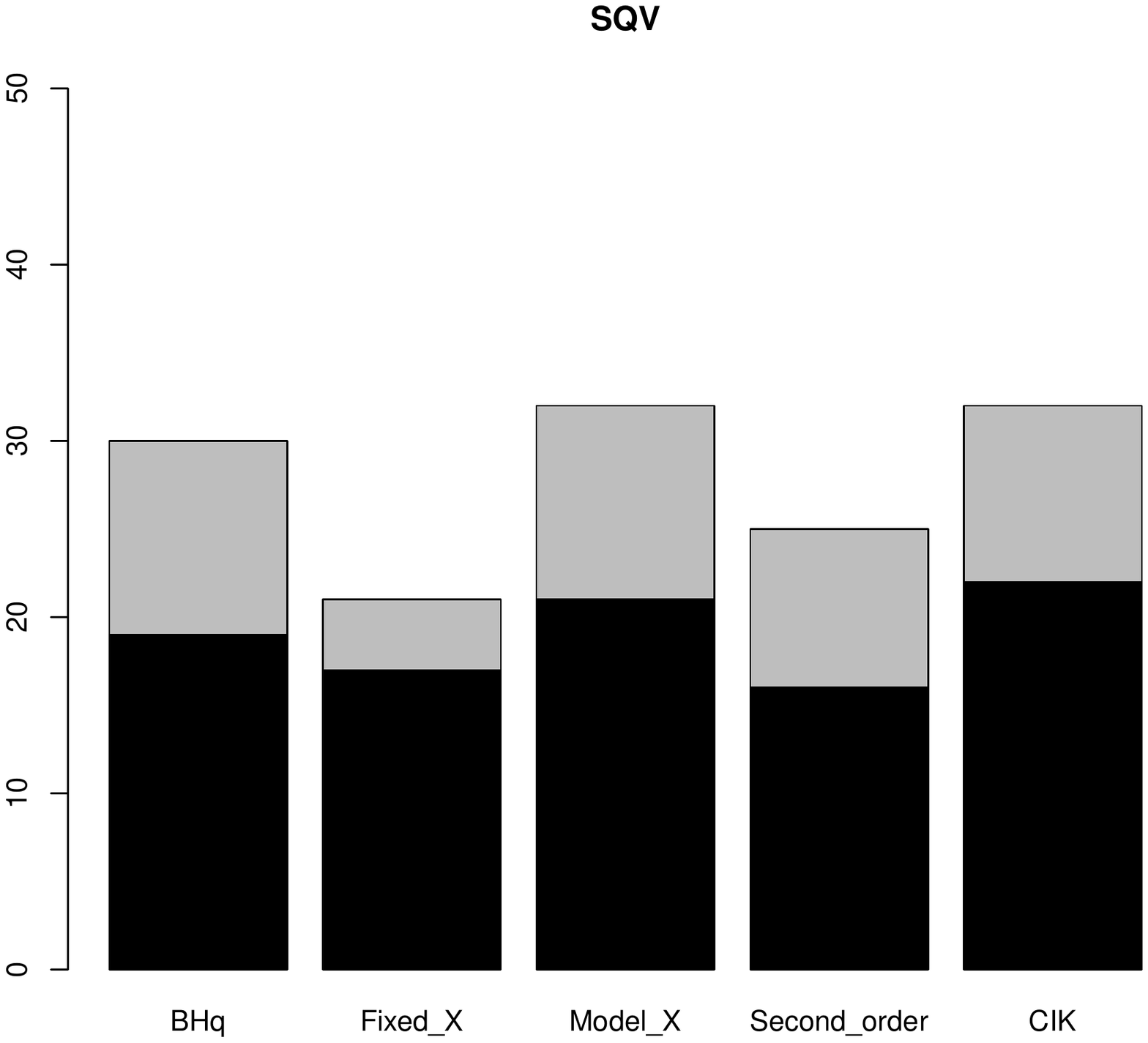}
           \caption{Comparison of different knockoff filters in terms of discoveries for each type of drug.  The black part denotes the true discoveries whilst the gray part denotes the false discoveries}
        \label{fig:seven drugs}
\end{figure}

\noindent Figure \ref{fig:real} summarizes the performances of the five methods across different drugs in terms of power and false discovery rate. It turns out that, in most cases, the CIKs are performing well. Compared to the other procedures, the CIKs are performing better in terms of power for APV, IDV and LPV whilst are performing worse for SQV. In terms of false discovery rate, the CIKs perform better than others for RTV whilst are performing worse for LPV, NFV and SQV. Figure \ref{fig:seven drugs} shows the performances of the five methods for each drug related to their discoveries. We note that the number of true discoveries with the CIKs is higher compared to BHq and Fixed-X for all the drugs and similarly to Second-order and Model-X. We also highlight the performance of the CIKs in RTV with respect to the other methods.

\medskip

\noindent To sum up, though the CIKs are not the best, they guarantee a good balance between power and false discovery rate and its performance is analogous to that of the other methods. For instance, as regards APV, ATV, IDV, LPV, NFV, the CIKs have a similar number of true discoveries with respect to Second-order and X-Model but also a fewer number of false discoveries.

\newpage

\bigskip

\noindent \textbf{Acknowledgments}
We are grateful to Guido Consonni for a very useful conversation.

%%%%%%%%%%%%%%%%%%%%%%%%%%%%%%%%%%%%%%%%%%%%%%
%% Supplementary Material, if any, should   %%
%% be provided in {supplement} environment  %%
%% with title and short description.        %%
%%%%%%%%%%%%%%%%%%%%%%%%%%%%%%%%%%%%%%%%%%%%%%

\noindent \textbf{Supplement}

\noindent We close the paper by proving Theorems \ref{t1} and \ref{t2}. For $c>0$ and $z\in\mathbb{R}^p$, we denote by $\phi_z(\cdot)$ the density function of $\mathcal{N}_p(z,\,c I)$, i.e.
$$\phi_z(x)=(2\,\pi\,c)^{-p/2}\exp\Bigl\{-\frac{\norm{x-z}^2}{2c}\Bigr\}\quad\quad\text{for all }x\in\mathbb{R}^p.$$

\medskip

\begin{proof}[Proof of Theorem \ref{t1}]
Define $c=\epsilon^2/2p$ and
$$\mu_0(A)=\int_{\mathbb{R}^p}\mathcal{N}_p(x,\,c I)(A)\,\mu(dx)\quad\quad\text{for all Borel sets }A\subset\mathbb{R}^p.$$
To see that $\mu_0\in\mathcal{P}_0$, it suffices to let
$$\Theta=\mathbb{R}^p,\quad\gamma=\mu\quad\text{and}\quad P_i(\cdot\mid\theta)=\mathcal{N}_1(\theta_i,c),$$
where $\theta_i$ denotes the $i$-th coordinate of $\theta=(\theta_1,\ldots,\theta_p)\in\mathbb{R}^p$. Obviously, $P_i(\cdot\mid\theta)$ is absolutely continuous. Moreover, since
$$\mathcal{N}_p(x,\,c I)=\mathcal{N}_1(x_1,\,c)\times\ldots\times\mathcal{N}_1(x_p,\,c)\quad\quad\text{for all }x=(x_1,\ldots,x_p)\in\mathbb{R}^p,$$
one obtains
$$\mu_0(A_1\times\ldots\times A_p)=\int_\Theta\prod_{i=1}^p\,P_i(A_i\mid\theta)\,\gamma(d\theta)\quad\quad\text{for all }A_1,\ldots,A_p\in\mathcal{B}.$$

\medskip

\noindent We next prove $d_{BL}(\mu,\mu_0)<\epsilon$. Fix a 1-Lipschitz function $g:\mathbb{R}^p\rightarrow [-1,1]$. Then,
\begin{gather*}
\int_{\mathbb{R}^p} g\,d\mu_0=\int_{\mathbb{R}^p}\int_{\mathbb{R}^p} g(y)\,\phi_x(y)\,dy\,\mu(dx).
\end{gather*}
Since $g$ is 1-Lipschitz and
$$\int_{\mathbb{R}^p}\norm{y-x}^2\,\phi_x(y)\,dy=p\,c,$$
it follows that
\begin{gather*}
\Abs{\int_{\mathbb{R}^p} g\,d\mu_0-\int_{\mathbb{R}^p} g\,d\mu}=\Abs{\int_{\mathbb{R}^p}\int_{\mathbb{R}^p} \bigl\{g(y)-g(x)\bigr\}\,\phi_x(y)\,dy\,\mu(dx)}
\\\le\int_{\mathbb{R}^p}\int_{\mathbb{R}^p}\abs{g(y)-g(x)}\,\phi_x(y)\,dy\,\mu(dx)
\\\le\int_{\mathbb{R}^p}\int_{\mathbb{R}^p}\norm{y-x}\,\phi_x(y)\,dy\,\mu(dx)
\\\le\int_{\mathbb{R}^p}\sqrt{\int_{\mathbb{R}^p}\norm{y-x}^2\,\phi_x(y)\,dy}\,\,\mu(dx)=\sqrt{p\,c}=\frac{\epsilon}{\sqrt{2}}.
\end{gather*}
Therefore,
$$d_{BL}(\mu,\mu_0)=\sup_g\,\Abs{\int_{\mathbb{R}^p} g\,d\mu_0-\int_{\mathbb{R}^p} g\,d\mu}\le\frac{\epsilon}{\sqrt{2}}<\epsilon.$$
\end{proof}

\bigskip

\begin{proof}[Proof of Theorem \ref{t2}]
Let $m$ denote the Lebesgue measure on $\mathbb{R}^p$. Suppose $\mu$ is absolutely continuous and denote by $f$ a density of $\mu$ (with respect to $m$). Given $\epsilon>0$, there is a function $f_0$ on $\mathbb{R}^p$ such that:

\begin{itemize}

\item $f_0$ is a probability density (with respect to $m$);

\item $\int_{\mathbb{R}^p}\abs{f(x)-f_0(x)}\,dx<\epsilon$;

\item $f_0$ is of the form $f_0=\sum_{j=1}^ka_j\,1_{R_j}$, where $k$ is a positive integer, $a_j>0$ a constant, and $R_j$ a bounded rectangle, i.e.
$$R_j=I_{1j}\times\ldots\times I_{pj}$$
where $I_{ij}$ is a bounded interval of the real line for each $i=1,\ldots,p$;

\end{itemize}
see e.g. Theorem (2.41) of \cite[p. 69]{FOLLAND}. Define $\mu_0$ as the probability measure on $\mathbb{R}^p$ with density $f_0$. Since $\mu$ and $\mu_0$ are both absolutely continuous,
$$d_{TV}(\mu,\mu_0)=(1/2)\,\int_{\mathbb{R}^p}\abs{f(x)-f_0(x)}\,dx<\epsilon/2<\epsilon.$$
Moreover, $\mu_0$ can be written as
$$\mu_0=\sum_{j=1}^ka_j\,m(R_j)\,\Bigl(\mathcal{U}_{1j}\times\ldots\times\mathcal{U}_{pj}\Bigr)$$
where $\mathcal{U}_{ij}$ is the uniform distribution on the interval $I_{ij}$. Hence, letting
$$\Theta=\{1,\ldots,k\},\quad\gamma\{\theta\}=a_\theta\,m(R_\theta)\quad\text{and}\quad P_i(\cdot\mid\theta)=\mathcal{U}_{i\theta},$$
one obtains
$$\mu_0(A_1\times\ldots\times A_p)=\sum_{\theta=1}^ka_\theta\,m(R_\theta)\,\prod_{i=1}^p\mathcal{U}_{i\theta}(A_i)=\int_\Theta\prod_{i=1}^pP_i(A_i\mid\theta)\,\gamma(d\theta)$$
for all $A_1,\ldots,A_p\in\mathcal{B}$. Hence, $\mu_0\in\mathcal{P}_0$.

\medskip

\noindent This proves the first part of the Theorem. To prove the second part, suppose $f$ is Lipschitz and define $\mu_0$ by \eqref{x48uh6tqa3} with
$c=\frac{1}{4p}\,\Bigl(\frac{\epsilon}{b\,m(B)}\Bigr)^2$, where $b$ is a Lipschitz constant for $f$ and $B\subset\mathbb{R}^p$ a Borel set satisfying $\mu(B^c)<\epsilon/2$ and $0<m(B)<\infty$. Since $\mu_0\in\mathcal{P}_0$, as shown in the proof of Theorem \ref{t1}, we have only to prove that $d_{TV}(\mu,\mu_0)<\epsilon$. The density $f_0$ of $\mu_0$ can be written as
$$f_0(x)=\int_{\mathbb{R}^p}\phi_x(y)\,f(y)\,dy.$$
Therefore,
\begin{gather*}
d_{TV}(\mu,\mu_0)=\int_{\mathbb{R}^p}(f(x)-f_0(x))^+dx\le\int_{B^c}f(x)\,dx+\int_B\abs{f(x)-f_0(x)}\,dx
\\=\mu(B^c)+\int_B\,\Abs{\int_{\mathbb{R}^p}\bigl\{f(x)-f(y)\bigr\}\,\phi_x(y)\,dy\,}\,dx
\\\le\mu(B^c)+\int_B\int_{\mathbb{R}^p}\abs{f(x)-f(y)}\,\phi_x(y)\,dy\,dx
\\\le\mu(B^c)+b\,\int_B\int_{\mathbb{R}^p}\norm{y-x}\,\phi_x(y)\,dy\,dx
\\\le\mu(B^c)+b\,\int_B\,\sqrt{\int_{\mathbb{R}^p}\norm{y-x}^2\,\phi_x(y)\,dy}\,\,\,dx
\\=\mu(B^c)+b\,\int_B\,\sqrt{p\,c}\,\,dx=\mu(B^c)+b\,m(B)\,\sqrt{p\,c}=\mu(B^c)+(\epsilon/2)<\epsilon
\end{gather*}
where the last inequality is because $\mu(B^c)<\epsilon/2$. This concludes the proof.
\end{proof}

\end{document}